\documentclass[11pt,reqno]{amsart}
\usepackage{amsmath,amsthm,amssymb}
\usepackage{nccmath}
\usepackage{a4wide}
\usepackage[numbers]{natbib}
\usepackage[utf8]{inputenc}
\usepackage{latexsym}
\usepackage{amssymb}
\usepackage[mathscr]{eucal}
\usepackage{epsfig}
\usepackage{lscape}
\usepackage{amssymb}
\usepackage{amsmath}
\usepackage{amsfonts}
\usepackage{graphicx}
\usepackage{xcolor}
\usepackage{yhmath}
\usepackage{mathdots}
\usepackage{MnSymbol}
\usepackage{epsfig}
\usepackage{lineno}
\usepackage{tikz}
\usepackage{bbold}
\usepackage{cite}
\usepackage[absolute]{textpos}

\newcommand{\NN}{\mathbb{N}}

\newcommand{\RR}{\mathbb{R}}

\newcommand{\supp}{\mathrm{supp}}

\newcommand{\mdeg}[0]{\mathrm{mdeg}}
\newcommand{\prodesc}[2]{\langle #1 , #2 \rangle}
\newlength{\wdth}

\newtheorem{theorem}{Theorem}[section]
\newtheorem{lemma}[theorem]{Lemma}
\newtheorem{corollary}[theorem]{Corollary}

\newtheorem{propo}[theorem]{Proposition}
\newtheorem{ex}[theorem]{Example}
\newtheorem{definition}[theorem]{Definition}
\newtheorem{remark}[theorem]{Remark}
\usepackage{mathtools}  
\usepackage{mathrsfs}   

\def\ps@pprintTitle{%
     \let\@oddhead\@empty
     \let\@evenhead\@empty
     \let\@evenfoot\@oddfoot}

\title{Multiple Orthogonal Polynomials of two real variables}

\author{Lidia Fern\'andez and Juan Antonio Villegas}

\address{Lidia Fern\'andez\\
IMAG and Departamento de Matem\'atica Aplicada\\
Universidad de Granada\\
18071, Granada, Spain.}
\email{lidiafr@ugr.es}
\address{Juan Antonio Villegas\\
IMAG and Departamento de Matem\'atica Aplicada\\
Universidad de Granada\\
18071, Granada, Spain}
\email{jantoniovr@ugr.es}
\thanks{The work of the authors was partially supported by IMAG-Mar\'ia de Maeztu grant CEX2020-001105-M.
}

\date{\today}



\begin{document}

\maketitle

\begin{abstract}
    Polynomials known as Multiple Orthogonal Polynomials in a single variable are polynomials that satisfy orthogonality conditions concerning multiple measures and play a significant role in several applications such as Hermite-Padé approximation, random matrix theory or integrable systems. However, this theory has only been studied in the univariate case. We give a generalization of Multiple Orthogonal Polynomials for two variables. Moreover, an extended version of some of the main properties are given. Additionally, some examples are given along the paper. 
\end{abstract}


\section{Introduction}\label{secIntro}

Multiple Orthogonality is a theory that extends standard orthogonality. In this framework, polynomials (defined on the real line in this introduction) satisfy orthogonality relations with respect to more than one measure. There are some results for Multiple Orthogonal Polynomials (MOP) that generalize their standard antecedents. For instance, the well-known theorem of existence and unicity of orthogonal polynomials sequences for a fixed real measure $\mu$ concerning its moments Hankel matrix has a generalization for Multiple Orthogonal Polynomials. Additionally, more extended results such as ``Nearest Neighbour Relations'' \citep{vanassche2011}, which generalize the three-term recurrence relation, will be shown in Section \ref{secPrelim}. 

Moreover, this extension has been widely studied recently, due to its applications in other mathematics branches \citep{martinezfinkelshtein-2016}. An introduction to multiple orthogonality and its properties can be found in \citep{Aptekarev}, \citep{Ismail} or \citep{VA20}. Some examples of applications of this theory are Hermite-Padé rational approximation \citep{vanassche2006}, Bessel functions and Bessel paths \citep{coussement-2003, Kuijlaars_2008}, electrostatics \citep{Andrei2022} or Markov chains and random walks \citep{Foulquie23, branquinho2021-2, branquinho2021}.

As shown, several researchers study and focus on this topic. However, Multiple Orthogonal Polynomials have been studied only in the univariate case. In this moment, there is no extension of the definitions and main results of multiple orthogonality employing bivariate polynomials, nor of its applications.

Nevertheless, it is possible to manage with standard orthogonal polynomials in the multivariate case. There are two different approaches. In 1975, T. Koornwinder employed the graded reverse lexicographic order in the usual basis of the space of multivariate polynomials to extend the standard definitions \citep{koornwinder}. This order has been employed to extend some families of multivariate orthogonal polynomials \citep{BP23}, and a similar order was applied to the Laurent monomials to study bivariate orthogonal Laurent polynomials \citep{CF24}. On the other hand, it is possible to bring polynomials of equal degree together and use polynomial vectors \citep{dunkl_xu_2014}. This approach has been employed to study a bi-dimensional version of Toda equations \citep{BP17}, or Quasi-birth-and-death processes \citep{FI21, FI23}, among others. More of these generalizations will be discused in Section \ref{SecBVMOP}. Some applications of orthogonal polynomials like Gaussian cubature or Padé approximation has also been generalized to the bivariate case \citep{CB00, Cuy83, Cuy86, Cuy99,CLY16}.

The structure of this article is the following. Firstly, in Section \ref{secPrelim}, we fix the notations and present the main definitions of multiple orthogonality in one variable: Type I and Type II MOP, together with a result which characterizes the existence of these polynomials. Next, in Section \ref{SecBVMOP}, two definitions of Type I and II bivariate MOP are shown. Additionally, these definitions are presented by using polynomial vectors. Following, some main results of the univariate case generalized to the bivariate one are presented in Section \ref{SecResults}: Existence and unicity of bivariate MOP, Biorthogonality relation and Nearest Neighbour Recurrence Relations. Previous to the conclusions, in Section \ref{SecRelation}, we find a relationship between the univariate and the bivariate cases, discussing the product of two univariate Type II MOP as a bivariate MOP.

\section{Preliminaries: Multiple Orthogonal Polynomials of one real variable}\label{secPrelim}

First, let us consider $r$ different real measures $\mu_1,\dots,\mu_r$ such that $\Omega_j=\supp(\mu_j)\subseteq\RR$, $j=1,\dots,r$. We will use multi-indices $\vec n = (n_1, \dots,n_r)\in \NN^r$, and denote $|\vec n| := n_1 + \dots + n_r$, the so called \textit{modulus} or \textit{norm} of $\vec n$. These multi-indices will determine the orthogonality relations with each measure and the degree of the polynomials.

For a fixed measure $\mu_j$, $j=1,\dots,r$, we denote its associated integral inner product as
\begin{equation}
    \label{eq:inner-product}
    \prodesc{f}{g}_j:=\int_{\Omega_j}f(x)g(x)d\mu_j(x).
\end{equation}

There are two different types of multiple orthogonality, which will be explained shortly. More information and a deeper introduction to multiple orthogonality can be found in \linebreak  \citep[Section 23.1]{Ismail} or \citep[Section 1]{VA20}. 
 
\begin{definition}[Type II Multiple Orthogonal Polynomials]
    \label{def:typeII-univar}
    Let $\vec n = (n_1,\dots,n_r)$. A monic polynomial $P_{\vec n}(x)$ is a \textbf{Type II Multiple Orthogonal Polynomial} if $\deg(P_{\vec n})= |\vec n|$ and 
    \begin{equation}
        \label{eq:typeII-MOP}
        \prodesc{P_{\vec n}}{x^k}_j = 0, \ \ \ k=0,\dots,n_{j}-1, \ \ j = 1,\dots,r
    \end{equation}
\end{definition}

This means $P_{\vec n}$ is orthogonal to $1,x,x^2,\dots,x^{n_j-1}$ with respect to each measure $\mu_j$, \linebreak $j=1,\dots,r$. 

On the other hand, there exists another type of multiple orthogonality: Type I multiple orthogonality.

\begin{definition}[Type I Multiple Orthogonal Polynomials]
    \label{def:typeI-univar}
    Let $\vec n = (n_1,\dots,n_r)$. \textbf{Type I Multiple Orthogonal Polynomials} are $r$ polynomials $A_{\vec n, 1}(x), \dots, A_{\vec n, r}(x)$, where  \linebreak $\deg(A_{\vec n, j})\leq n_j-1$, $j=1,\dots,r$ and satisfying
    \begin{equation}
        \label{eq:typeI-MOP}
        \sum_{j=1}^r \prodesc{A_{\vec n,j}}{x^k}_j = \left\{\begin{array}{ccl}
            0 &   \text{ if } & k=0,\dots,|\vec n|-2 \\
            1 & \text{ if } & k=|\vec n|-1.    
        \end{array}\right.
    \end{equation}
    
\end{definition}

Whenever the measures are all absolutely continuous with respect to a common positive measure $\mu$ defined in $\Omega = \displaystyle\bigcup_{i=1}^r \Omega_i$, \textit{i.e.}, $d\mu_j = w_j(x) d\mu(x)$, $j=1,\dots,r$, it is possible to define the \textit{Type I function} as
\begin{equation}
    \label{eq:typeI-function}
    Q_{\vec n}(x)=\sum_{j=1}^r A_{\vec n,j}(x)w_j(x).
\end{equation}
Using the Type I function, the orthogonality relations \eqref{eq:typeI-MOP} could be written as
\begin{equation}
    \label{eq:typeI-MOP-function}
    \prodesc{Q_{\vec n}}{x^k}_{\mu}= \left\{\begin{array}{ccl}
        0 &   \text{ if } & k=0,\dots,|\vec n|-2 \\
        1 & \text{ if } & k=|\vec n|-1.      
    \end{array}\right.
\end{equation}

Nevertheless, neither all systems of measures $\mu_1,\dots,\mu_r$ nor all multi-indices $\vec n\in\NN^r$ yield a Type I vector of polynomials or a Type II polynomial, just as in standard OP theory, where not every measure produces an OPS $\{P_n\}_{n\geq 0}$, depending on the Hankel matrix of moments. In this sense, there is an analogous result for existence and unicity of MOP.
\begin{propo}\textnormal{\citep[Section 23.1]{Ismail}}
    \label{prop:existence-of-MOP}
    Given a multi-index $\vec n\in\NN^r$ and $r$ positive measures, $\mu_1,\dots,\mu_r$, the following statements are equivalent:
    \begin{enumerate}
        \item There exist a unique vector $(A_{\vec n,1}, \dots, A_{\vec n,r})$ of Type I MOP.
        \item There exist a unique Type II MOP $P_{\vec n}$.
        \item The following block matrix of moments is regular: 
        \begin{equation}
            \label{eq:MOP-matrix}
            M_{\vec n}=\left(\begin{array}{c}
                M_{n_1}^{(1)} \\ \hline
                M_{n_2}^{(2)} \\ \hline
                \vdots \\ \hline
                M_{n_r}^{(r)} \\ 
            \end{array}\right), \text{ \ \  where \ \ } M_{n_j}^{(j)} = \begin{pmatrix}
                m_0^{(j)} & m_1^{(j)} & \cdots & m_{|\vec n|-1}^{(j)} \\
                m_1^{(j)} & m_2^{(j)} & \cdots & m_{|\vec n|}^{(j)} \\
                \vdots & \vdots & \ddots & \vdots \\
                m_{n_j-1}^{(j)} & m_{n_j}^{(j)} & \cdots & m_{|\vec n|+n_j-2}^{(j)} \\
            \end{pmatrix}, \ \ j=1,\dots,r, 
        \end{equation}
        and $m_k^{(j)}=\displaystyle\int_{\Omega_j} x^k d\mu_j(x)$ are the moments of the measures, $k\geq 0$.
    \end{enumerate}
\end{propo}

This means that, for a fixed multi-index, Type I MOP associated to this multi-index exists if and only if Type II MOP exists. Following Proposition \ref{prop:existence-of-MOP}, we provide the next definition.

\begin{definition}
    A multi-index $\vec n = (n_1,\dots,n_r)\in\NN^r$ is \textbf{normal} if it satisfies the conditions in Proposition \ref{prop:existence-of-MOP}.
    A system of $r$ measures $\mu_1,\dots,\mu_r$ is \textbf{perfect} if every $\vec n\in\NN^r$ is normal.
\end{definition}

Some examples of perfect systems are Angelesco systems, AT-systems and Nikishin systems, see \citep[Sections 23.1.1 and 23.1.2]{Ismail}, \citep{nikishin1991rational}.

\section{Multiple Orthogonal polynomials of two real variables}\label{SecBVMOP}

In the bivariate case, it is possible to deal with Orthogonal Polynomials by using two different approaches. In \citep{dunkl_xu_2014}, C. Dunkl and Y. Xu employ a polynomial vector notation, so that, if we denote as $\mathbb X_i = (x^i, x^{i-1}y, \dots, y^i)^t, i\geq 0,$ the vector of all monomials of total degree $i$, whose size is $i+1$, we can express a polynomial vector of degree $n$ as
$$
\mathbb P_n = \sum_{i=0}^n G_{i,n} \mathbb X_i,
$$
where $G_{i,n}$ is a $(n+1)\times (i+1)$ matrix and $G_{n,n}$ is regular. With this notation, it is \linebreak  possible to extend the standard orthogonality definitions and results to polynomial vectors, \linebreak  see \citep[Chapter 3]{dunkl_xu_2014}.

A different approach can be found in \citep{koornwinder}, where T. Koornwinder uses the basis of the vector space of bivariate polynomials ordered by the graded reverse lexicographic order, this is
\begin{equation}
    \label{eq:usual-basis-ordered}
    \{1,x,y,x^2,xy,y^2,x^3,x^2y,xy^2,y^3,\dots\}.
\end{equation}
Our extension to the bivariate case of multiple orthogonality is mainly focused in this second approach, although some properties will be also presented by using polynomial vectors in order to clarify some notations.

\begin{figure}[h]
    \centering\includegraphics[width=5cm]{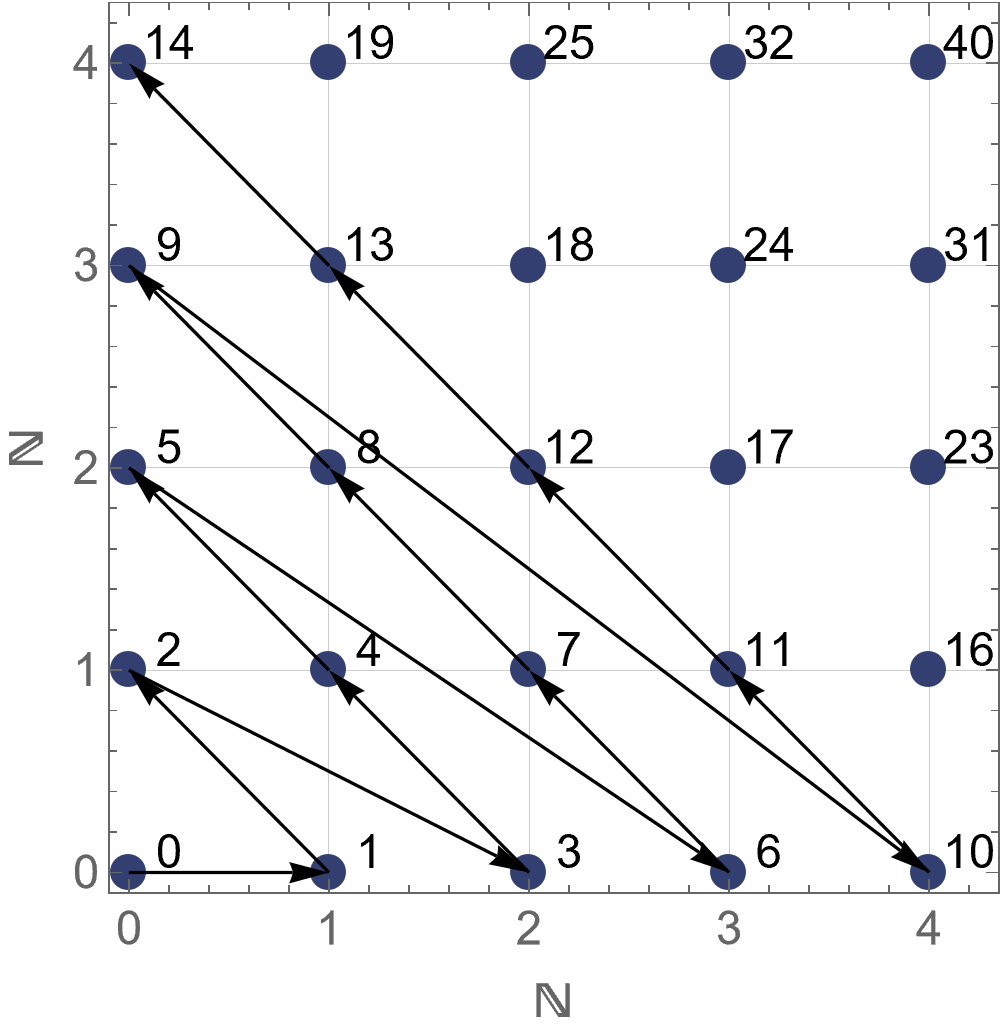}
    \caption{Graphic representation of Cantor Pairing Function}
    \label{fig:Cantor}
\end{figure}

Previous to our definitions, remember that the \emph{Cantor Pairing Function}, $\pi$, is a bijection between $\NN$ and $\NN^2$ (see \citep{Lisi}). Its expression is given by 
\begin{equation}
    \label{eq:CPF}
\pi(t, s) = \dfrac{1}{2}(t + s)(t + s + 1) + s.
\end{equation}
In Figure  \ref{fig:Cantor} you can observe a graphic representation of this bijection.

Using the Cantor Pairing Function, it is possible to assign a natural number to each monomial of the basis \eqref{eq:usual-basis-ordered}, so that $1$ is associated to $0$, $x$ to $1$, $y$ to $2$, and so on. Following this idea, given a bivariate polynomial $P(x,y)$, you can express $P$ as a linear combination of the first elements of \eqref{eq:usual-basis-ordered}, in a way that  
$$P(x,y)= c_{\pi(l,m)} x^l y^m + c_{\pi(l,m)-1} x^{l+1} y^{m-1} + \cdots + c_2 y + c_1 x + c_0, \ \ \ c_{\pi(l,m)}\not=0.$$ 
We will refer as the \emph{multidegree} of $P(x,y)$, denoted by $\mdeg(P(x,y))$, to the pair $(l,m)$, which are the exponents of the variables `$x$' and `$y$' (in that order) of the leading term. Applying the function $\pi$ to that pair, observe $\pi(\mdeg(P))=\pi(l,m)$ is the position of the leading monomial $x^l y^m$ in the list \eqref{eq:usual-basis-ordered}, starting the count at $0$. Next, we present some useful definitions associated to every multi-index.

\begin{definition}[Degree and Remainder]
    \label{def:dnkn}
For every multi-index $\vec n =(n_1, \dots, n_r)\in\NN^r$, let us consider $|\vec n|$, the modulus of $\vec n$, and let $(l,m)$ be the pair such that $\pi(l,m)=|\vec n|$. We define the \textbf{degree} $d_{\vec n}$ and the \textbf{remainder} $k_{\vec n}$ associated to $\vec n$ as
\begin{equation}
    \begin{aligned}
    d_{\vec n} &= l+m, & k_{\vec n}&=|\vec n|-\frac 1 2 d_{\vec n}(d_{\vec n}+1)=m.
    \end{aligned}
\end{equation}
This means that it is possible to write
$$
|\vec n| = \frac 1 2 d_{\vec n}(d_{\vec n}+1)+k_{\vec n}.
$$
\end{definition}

\begin{remark} For every $l,m\in\NN$ and every multi-index $\vec n\in\NN^r$ such that $|\vec n|=\pi(l,m)$
    $$
    \pi(l,m) =  \pi(d_{\vec n}-k_{\vec n},k_{\vec n})
    $$
\end{remark}

\begin{ex}
    In Table \ref{tab:parameters} some examples of multi-indices $\vec n\in\NN^r$ for some values of $r$ and their associated values of the modulus, multidegree, degree and remainder are shown. 
    \begin{table}[h]
        \begin{tabular}{cccccc}
        \hline
        $r$ & $\vec n$      & $|\vec n|$ & $\pi^{-1}(|\vec n|)$ & $d_{\vec n}$ & $k_{\vec n}$ \\ \hline
        $2$ & $(6,2)$       & $8$        & $(1,2)$              & $3$          & $2$          \\ \hline
        $3$ & $(2,1,1)$     & $4$        & $(1,1)$              & $2$          & $1$          \\ \hline
        $4$ & $(4,6,7,3)$   & $20$       & $(0,5)$              & $5$          & $5$          \\ \hline
        $5$ & $(1,6,2,1,2)$ & $12$       & $(2,2)$              & $4$          & $2$          \\ \hline
        \end{tabular}
        \vspace{0.3cm}
        \caption{Parameters associated to multi-indices}
        \label{tab:parameters}
        \end{table}
\end{ex}

Remember that, in the univariate case, $|\vec n|$ was the number of orthogonality conditions of both Type I and II MOP, and it was also the degree of the Type II MOP. In our generalization, given below, $|\vec n|$ is also the number of orthogonality conditions, but, instead of the degree of the bivariate Type II MOP, in this case it is the position of its leading term in \eqref{eq:usual-basis-ordered}, \textit{i.e.}, $|\vec n| = \pi(\mdeg(P_{\vec n}))$. As a consequence, $d_{\vec n}$ is the actual degree of $P_{\vec n}$.

Now, given $r\in\NN$, and a system of $r$ $2$-dimensional measures $\mu_1, \dots, \mu_r$, we present the definitions of Type I and II multiple orthogonality.

\begin{definition}[Type II Bivariate Multiple Orthogonality]
    \label{def:typeII}
    Let $d\in\NN$ and $l,m\in\NN$ such that $l+m=d$. Take a multi-index $\vec n=(n_1,\dots,n_r)\in\NN^r$ such that $|\vec n|=\pi(l,m)$. Then, the monic polynomial  
    $$
    P_{\vec n}(x,y) = x^l y^m + c_{|\vec n|-1} x^{l+1} y^{m-1} + \cdots + c_1 x + c_0,
    $$
    is the \textbf{bivariate Type II Multiple Orthogonal Polynomial} associated to $\vec n$ if it satisfies the orthogonality conditions
  \begin{equation}
    \label{eq:typeII-MOP-2-app}
    \prodesc{P_{\vec n}}{x^t y^s}_j = 0, \ \ \ 0\leq \pi(t,s)\leq n_j-1, \ j=1,\dots,r.
  \end{equation}    
\end{definition}

Observe that, as mencioned, $\pi(\mdeg(P_{\vec n})) = |\vec n|$ and $\deg(P_{\vec n})=l+m=d_{\vec n}$, where $d_{\vec n}$ is the degree associated to $\vec n$ (see Definition \ref{def:dnkn}).

\begin{definition}[Type I Bivariate Multiple Orthogonality]
    \label{def:typeI}
    Let $\vec n=(n_1,\dots,n_r)\in\NN^r$ be a multi-index. For each $n_j$, consider $l_j,m_j\in\NN$ such that $\pi(l_j,m_j)=n_j-1$, $j=1,\dots,r$. Then, the polynomials $A_{\vec n,1}(x,y),\dots,A_{\vec n,r}(x,y)$ given by
    $$A_{\vec n,j} = c_{n_j-1,j} x^{l_j}y^{m_j}+ c_{n_j-2,j}x^{l_j+1}y^{m_j-1} +\cdots +c_{1,j} x + c_{0,j}, \ \ \ j=1,\dots,r,$$
    are the \textbf{bivariate Type I Multiple Orthogonal Polynomials} if they satisfy the orthogonality conditions
      \begin{equation}
        \label{eq:condition-type-I-2-app}
        \sum_{j=1}^r \prodesc{A_{\vec n,j}}{x^t y^s}_j = \left\{\begin{array}{ccl}
            0 &   \text{ if } & 0\leq \pi(t,s)\leq |\vec n|-2 \\
            1 & \text{ if } & \pi(t,s)= |\vec n|-1     
        \end{array}\right..
      \end{equation}
\end{definition}

In this case, we have $\pi(\mdeg(A_{\vec n,j}))\leq n_j-1$, $j=1,\dots,r$. 

Once again, let us assume the measures $\mu_1,\dots,\mu_r$ are absolutely continuous with respect to a common measure $\mu$ whose support is $\Omega = \displaystyle\bigcup_{j=1}^r \Omega_j$, where $\Omega_j=\supp(\mu_j)$ (this is, \linebreak  $d\mu_j(x,y)=w_j(x,y)d\mu(x,y), \ j=1,\dots,r$). Then, we define the \textit{Bivariate Type I Function} as
\begin{equation}
    \label{eq:typeI-function-2-app}
    Q_{\vec n}(x,y)=\sum_{j=1}^r A_{\vec n,j}(x,y) w_j(x,y),
\end{equation}
so that conditions \eqref{eq:condition-type-I-2-app} are equivalent to
\begin{equation}
    \label{eq:typeI-conditions-function-2-app}
    \prodesc{Q_{\vec n}}{x^t y^s}_\mu = \left\{\begin{array}{ccl}
        0 &   \text{ if } & 0\leq \pi(t,s)\leq |\vec n|-2 \\
        1 & \text{ if } & \pi(t,s)= |\vec n|-1     
    \end{array}\right..
\end{equation}

In order to clarify this definitions, we will introduce the following example.

\begin{ex}
    Consider a system of $r=3$ bidimensional measures. If we choose $d=3, l=2$ and $m=1$, we have to choose a multi-index $\vec n\in\NN^3$ such that $|\vec n|=\pi(2,1)=7$, for example $\vec n =(2,4,1)$. As $|\vec n| = 7 = \frac 1 2 \cdot 3\cdot 4 + 1$, then $d_{(2,4,1)} = 3$ and $k_{(2,4,1)} = 1$.
    Assuming the existence, which will be discussed later, the bivariate Type II MOP will have the following shape:
    $$
    P_{(2,4,1)}(x,y)=x^2 y + c_6 x^3 + c_5 y^2 + c_4 x y + c_3 x^2 + c_2 y + c_1 x + c_0,
    $$
    and it will satisfy
    $$\prodesc{P_{(2,4,1)}}{1}_1 = \prodesc{P_{(2,4,1)}}{x}_1 = 0,$$
    $$\prodesc{P_{(2,4,1)}}{1}_2 = \prodesc{P_{(2,4,1)}}{x}_2 =\prodesc{P_{(2,4,1)}}{y}_2 =\prodesc{P_{(2,4,1)}}{x^2}_2 = 0,$$
    $$\prodesc{P_{(2,4,1)}}{1}_3 = 0.$$

    Note that, indeed, $\deg(P_{(2,4,1)}) = d_{(2,4,1)} = 3$.

    On the other hand, the bivariate Type I MOPs will have the following shape:
    \begin{equation*}
        \begin{array}{c}
            A_{(2,4,1),1}(x,y)= c_{1,1} x + c_{0,1} \\
            A_{(2,4,1),2}(x,y)= c_{3,2} x^2 + c_{2,2} y + c_{1,2} x + c_{0,2} \\
            A_{(2,4,1),3}(x,y)= c_{0,3} \\
        \end{array}
    \end{equation*}
    and they will satisfy
    \begin{equation*}
        \sum_{j=1}^3 \prodesc{A_{\vec n,j}}{1}_j =  \sum_{j=1}^3 \prodesc{A_{\vec n,j}}{x}_j = \sum_{j=1}^3 \prodesc{A_{\vec n,j}}{y}_j = \sum_{j=1}^3 \prodesc{A_{\vec n,j}}{x^2}_j = \sum_{j=1}^3 \prodesc{A_{\vec n,j}}{xy}_j = \sum_{j=1}^3 \prodesc{A_{\vec n,j}}{y^2}_j = 0,
    \end{equation*}
      $$
      \sum_{j=1}^3 \prodesc{A_{\vec n,j}}{x^3}_j = 1
      $$

\end{ex}

In order to bring polynomials of equal degree together and clarify some notations according to \citep{dunkl_xu_2014}, we extended Definitions \ref{def:typeII} and \ref{def:typeI} to polynomial vectors, which we will call bivariate Multiple Orthogonal Polynomial Vectors (MOPV). But first, let us establish some conditions regarding the multi-indices and parameters that will be employed. Given $d\in\NN$, we know there are $d+1$ different monomials of degree $d$. Following the order established by the application $\pi$, the first multidegree is $(d,0)$ (associated to the monomial $x^d$). From this point, choose $d+1$ multi-indices $\vec n_0,\dots,\vec n_d$ such that
    \begin{itemize}

        \item The modulus of the multi-indices satisfy
        $$
        \begin{array}{l}
            |\vec n_0| = \pi(d,0) = \frac 1 2 d(d+1),\\[3pt]
            |\vec n_1| = \pi(d-1,1) = \frac 1 2 d(d+1)+1\\[3pt]
            \qquad \vdots \\[3pt]
            |\vec n_d| = \pi(0,d) = \frac 1 2 d(d+1)+d\\[3pt]              
        \end{array}
        $$
        This is,
        \begin{equation}
            \label{eq:condition-multi-indices-1}
            |\vec n_k|=\frac 1 2 d(d+1)+k, \ \  k=0,\dots,d.
        \end{equation}
        As a result, every multi-index is associated to a bivariate Type II MOP of degree $d$.
        \item They are neighbour multi-indices, \textit{i.e.}
        \begin{equation}
            \label{eq:condition-multi-indices-2}
            \vec n_{k+1} = \vec n_k + \vec e_j, \ \ k=0,\dots,d-1,
        \end{equation}
        where $\vec e_j = (0,\dots,\underset{(j)}{1},\dots,0)$, for certain $j\in\{1,\dots,r\}$.
    \end{itemize}

    Under these conditions, it is possible to construct a degree $d$ polynomial vector, employing $d+1$ multi-indices satisfying conditions \eqref{eq:condition-multi-indices-1} and \eqref{eq:condition-multi-indices-2}. We will denote as $\vec n_k =(n_{k,1},\dots,n_{k,r})$ the components of the multi-index $\vec n_k$, $k=0,\dots,d$.

    \begin{definition}[Bivariate Type II Multiple Orthogonal Polynomial Vector]
    \label{def:typeII-vector}
    Given $d\in\NN$, consider $d+1$ multi-indices  $\vec n_0,\dots,\vec n_d$  satisfying conditions \eqref{eq:condition-multi-indices-1} and \eqref{eq:condition-multi-indices-2}. For each multi-index $\vec n_k$, let $P_{\vec n_k}$ be its associated bivariate Type II MOP of Definition \ref{def:typeII}. The polynomial vector 
    $$
    \mathbb{P}_{\{\vec n_k\}_0^d}^{(d)} = \begin{pmatrix}
      P_{\vec n_0} \\ \vdots \\ P_{\vec n_d}\end{pmatrix} = \sum_{k=0}^{d} G_{d,k} \mathbb X_k,
    $$
    where $G_{d,k}$ is a $(d+1)\times(k+1)$ matrix and $G_{d,d}$ is a lower triangular matrix with ones on the diagonal, is the \textbf{Bivariate Type II Multiple Orthogonal Polynomial Vector} associated to $\{\vec n_0,\dots,\vec n_d\}$ and it satisfies 
       
		\begin{equation*}
		  \langle\mathbb{P}_{\{\vec n_k\}_0^d}^{(d)},(\mathbb X_0^t | \cdots | \mathbb X_d^t)\rangle_j = \left(\begin{array}{c}
            R_0 \\ \hline \vdots \\  \hline R_d 
          \end{array}\right),
        \end{equation*}
        where the $k$-th row $R_k$ has its $n_{k,j}$ first components equal to $0$, $k=0,\dots,d$, $j=1,\dots,r$.
    \end{definition}

    \begin{remark}
        The polynomial vector $\mathbb{P}_{\{\vec n_k\}_0^d}^{(d)}$ is composed of polynomials of degree $d$ so that their leading terms are $x^d, x^{d-1}y,\dots,y^d$.
    \end{remark}

    In order to clarify this definition, we present an example.
    \begin{ex}
        \label{example:typeIIMOPV}
        Let us take $r=2$, $d=2$, then, we need to choose three neighbour multi-indices $\vec n_0, \vec n_1, \vec n_2$ so that $|\vec n_0|=\pi(2,0)=2\cdot 3/2 + 0 = 3, |\vec n_1|=\pi(1,1)=4, |\vec n_2|=\pi(0,2)=5$. Consider, for example, $\vec n_0=(1,2), \vec n_1 =(1,3), \vec n_2 = (2,3)$. Then
	$$ \mathbb P_{\{\vec n_0, \vec n_1, \vec n_2\}}^{(2)} = (P_{(1,2)},P_{(1,3)},P_{(2,3)})^t$$
	satisfying
	$$ \prodesc{\mathbb P_{\{\vec n_0, \vec n_1, \vec n_2\}}^{(2)}}{(1,x,y,x^2,xy,y^2)}_1 = \begin{pmatrix}
		0 & \ast & \ast& \ast& \ast& \ast \\
		0 & \ast & \ast& \ast& \ast& \ast \\
		0 & 0 & \ast& \ast& \ast& \ast \\ 
	  \end{pmatrix}, $$
	since the first components of each pair are 1, 1 and 2.
    
      $$
	\prodesc{\mathbb P_{\{\vec n_0, \vec n_1, \vec n_2\}}^{(2)}}{(1,x,y,x^2,xy,y^2)}_2 = \begin{pmatrix}
		0 & 0 & \ast& \ast& \ast& \ast \\
		0 & 0 & 0 & \ast& \ast& \ast \\
		0 & 0 & 0 & \ast& \ast& \ast \\ 
	  \end{pmatrix},
	$$
    since the second components of each pair are 2, 3 and 3.
    \end{ex}

    Following this idea, we could also define the bivariate Type I MOPVs. In this case, we define $r$ polynomial vectors of different degree but necessarily same size. 

    \begin{definition}[Type I Multiple Orthogonal Polynomial Vectors] 
        \label{def:typeI-vector}
        Given $d\in\NN$, consider $d+1$ multi-indices  $\vec n_0,\dots,\vec n_d$ satisfying conditions \eqref{eq:condition-multi-indices-1} and \eqref{eq:condition-multi-indices-2}. For each multi-index $\vec n_k$, let $A_{\vec n_k,1},\dots,A_{\vec n_k,r}$ be its associated bivariate Type I MOP of Definition \ref{def:typeI}. The polynomial vectors
        $$
        \mathbb{A}_{\{\vec n_k\}_0^d,j} = \begin{pmatrix}
          A_{\vec n_0,j} \\ \vdots \\ A_{\vec n_d,j}\end{pmatrix}, \ \ \ j=1,\dots,r
        $$
        are the \textbf{Bivariate Type I Multiple Orthogonal Polynomial Vectors} associated to $\{\vec n_0,\dots,\vec n_d\}$, and they satisfy 
        \begin{equation*}
            \sum_{j=1}^r\prodesc{\mathbb{A}_{\{\vec n_k\}_0^d,j}}{(\mathbb X_0^t | \cdots | \mathbb X_d^t)}_j = \left(\begin{array}{c}
                \tilde R_0 \\ \hline \vdots \\  \hline \tilde R_d  
              \end{array}\right),
        \end{equation*}
        where the $k$-th row $\tilde R_k$ has its $|\vec n_k|-1$ first components equal to $0$ and the $|\vec n_k|$-th component is $1$, $k=0,\dots,d$, $j=1,\dots,r$.
        
      \end{definition}

    \begin{remark}
        As mentioned, the polynomial $A_{\vec n_k, j}$ is the $j$-th Type I MOP of Definition \ref{def:typeI} for the multi-index $\vec n_k$. Then, $\pi(\mdeg(A_{\vec n_k,j})) \leq n_{k,j}-1$. Besides, observe that all polynomial vectors have the same size ($d+1$) and there are as many bivariate Type I MOPV as measures, $r$.
    \end{remark}
    Naturally, the \textit{Bivariate Type I Vector} is defined as 
    \begin{equation}
        \label{eq:TypeIVector}
        \mathbb{Q}_{\{\vec n_k\}_0^d} = \sum_{j=1}^r \mathbb{A}_{\{\vec n_k\}_0^d,j}w_j(x,y).
    \end{equation}

    Observe the following example, where we use the same parameters and the same multi-indices as in Example \ref{example:typeIIMOPV}.
    \begin{ex}
        Let us take $r=2$, $d=2$, then $|\vec n_0|= 3, |\vec n_1|=4, |\vec n_2|=5$, and we consider $\vec n_0=(1,2), \vec n_1 =(1,3), \vec n_2 = (2,3)$. Then, the bivariate Type I MOPV are

	\begin{equation*}
		\begin{aligned}
			\mathbb A_{\{\vec n_0, \vec n_1, \vec n_2\},1}&=\begin{pmatrix}
				A_{(1,2), 1} \\ A_{(1,3), 1} \\ A_{(2,3), 1} \\
			\end{pmatrix} & \mathbb A_{\{\vec n_0, \vec n_1, \vec n_2\},2}&=\begin{pmatrix}
				A_{(1,2), 2} \\ A_{(1,3), 2} \\ A_{(2,3), 2} \\
			\end{pmatrix}  
		\end{aligned}
	\end{equation*}
	and they satisfy
	$$ \sum_{j=1}^r \prodesc{\mathbb A_{\{\vec n_0, \vec n_1, \vec n_2\},j}}{(1,x,y,x^2,xy,y^2)}_j = \begin{pmatrix}
		0 & 0 & 1 & \ast& \ast& \ast \\
		0 & 0 & 0 & 1 & \ast& \ast \\
		0 & 0 & 0 & 0 & 1 & \ast \\ 
	  \end{pmatrix},
	$$
    since  $|\vec n_0|= 3, |\vec n_1|=4, |\vec n_2|=5$.
    \end{ex}

In the next sections, we will present generalized versions of some results from the univariate case employing Definitions \ref{def:typeII} and \ref{def:typeI}, and some analogous when the vector notation of Definitions \ref{def:typeII-vector} and \ref{def:typeI-vector} is applied.

\section{Main Results}\label{SecResults}

We will start with the existence of bivariate Type I and/or Type II MOP. Given a multi-index $\vec n=(n_1,\dots,n_r)$, we will denote by $M_{\vec n}$ the following block matrix:
\begin{equation}
    \label{eq:matrix-Mn}
    M_{\vec n} = \begin{pmatrix}
        M_{n_1}^{(1)} \ | \ \cdots \  | \  M_{n_r}^{(r)}
       \end{pmatrix},
\end{equation} 
where, for every $j=1,\dots,r$:
\begin{equation}
    \label{eq:matrix-Mn-blocks}
    M_{n_j}^{(j)} = \begin{pmatrix}
        m_{(0,0)}^{(j)} & m_{(1,0)}^{(j)} & \cdots & m_{\pi^{-1}(n_j-2)}^{(j)}  & m_{\pi^{-1}(n_j-1)}^{(j)}  \\
        m_{(1,0)}^{(j)} & m_{(2,0)}^{(j)} & \cdots & m_{\pi^{-1}(n_j-2)+(1,0)}^{(j)}  & m_{\pi^{-1}(n_j-1)+(1,0)}^{(j)}  \\
        \vdots & \vdots & \ddots & \vdots & \vdots \\
        m_{\pi^{-1}(|\vec n|-1)}^{(j)} & m_{\pi^{-1}(|\vec n|-1)+(1,0)}^{(j)} & \cdots & m_{\pi^{-1}(n_j-2)+\pi^{-1}(|\vec n|-1)}^{(j)}  & m_{\pi^{-1}(n_j-1)+\pi^{-1}(|\vec n|-1)}^{(j)}  \\
      \end{pmatrix},
\end{equation}
where $m_{(t,s)}^{(j)}=\displaystyle\int_{\mathbb R^2}x^t y^s d \mu_j(x,y)$ are the moments of the bidimensional measures. Observe that the value in the $k$-th row and the $l$-th column of $M_{n_j}^{(j)}$ (starting the count in 0) is $m_{\pi^{-1}(k)+\pi^{-1}(l)}^{(j)}$, $k=0,\dots,|\vec n|-1$, $l=0,\dots,n_j-1$.

For bivariate Type II MOP, applying conditions \eqref{eq:typeII-MOP-2-app} to a generic monic polynomial, it is easy to deduce that the coefficients of this polynomial $\mathbf{c}_{\vec n}=(c_0,\dots,c_{|\vec n|-1})^t$ are the solutions of the linear system
\begin{equation*}
    \label{eq:system-typeII-2app}
    M_{\vec n}^t\cdot \mathbf{c}_{\vec n} = \mathbf{b}_{\vec n}, \text{ where } \mathbf{b}_{\vec n} = - (\mathbf{b}_{\vec n}^{(1)} \ | \ \dots \ | \ \mathbf{b}_{\vec n}^{(r)})^t,
\end{equation*}
with 
\begin{equation*}
    \mathbf{b}_{\vec n}^{(j)} = \left(m_{\pi^{-1}(|\vec n|)}^{(j)},  m_{\pi^{-1}(|\vec n|)+(1,0)}^{(j)},  \dots , m_{\pi^{-1}(|\vec n|)+\pi^{-1}(n_j-1)}^{(j)}\right)  \ \ j=1,\dots,r.
\end{equation*}

Besides, applying \eqref{eq:condition-type-I-2-app} to generic bivariate Type I MOPs, we obtain that the coefficients of bivariate Type I MOPs $\mathbf{\tilde c}_{\vec n}= (\tilde c_{0,1},\dots,\tilde c_{n_1-1,1} | \ \dots \ |\tilde  c_{0,r},\dots,\tilde c_{n_r-1,r})^t$ are the solutions of the system 

\begin{equation*}
    \label{eq:system-typeI-2app}
    M_{\vec n} \cdot \mathbf{\tilde c}_{\vec n} = (0, \dots, 0, 1)^t.
\end{equation*}

Then, for a fixed multi-index $\vec n=(n_1,\dots,n_r)$, Type I and Type II will exist if and only if the matrix $M_{\vec n}$ is regular. So, we have the following result, which clearly generalizes Proposition~\ref{prop:existence-of-MOP}.
\begin{propo}
    \label{prop:equivalence-2-app}
    Given a system of $r$ bidimensional measures $\mu_1,\dots,\mu_r$ and a multi-index $\vec n=(n_1,\dots,n_r)\in\NN^r$, the following statements are equivalent.
    \begin{enumerate}
        \item There exist unique Type I polynomials $A_{\vec n,1},\dots,A_{\vec n,r}$.
        \item There exist a unique Type II polynomial $P_{\vec n}$.
        \item The matrix $M_{\vec n}$ defined in \eqref{eq:matrix-Mn} and \eqref{eq:matrix-Mn-blocks}  is regular.
    \end{enumerate}
\end{propo}

As happened in the univariate case, we have a matrix whose regularity for a fixed multi-index ensures the existence and unicity of both Type I and II bivariate MOP. Additionally, for a fixed multi-index, bivariate Type I MOPs exist if, and only if, bivariate Type II MOP exists. This proposition motivates the following definitions.

\begin{definition}
    Given a system of $r$ bidimensional measures $\mu_1,\dots \mu_r$ and a multi-index $\vec n=(n_1,\dots,n_r)\in\NN^r$, $\vec n$ is said to be \textbf{normal} if it satisfies the conditions in Proposition \ref{prop:equivalence-2-app}, and the system of measures $\mu_1,\dots,\mu_r$ is a \textbf{perfect} system of measures if every multi-index $\vec n\in\NN^r$ is normal. 
\end{definition}

See the following example to clarify the structure of these matrices.

\begin{ex}
    Let us consider a system of $r=3$ bidimensional measures $\mu_1,\mu_2,\mu_3$ and a multi-index $\vec n =(1,2,2)$. Then, bivariate Type I and Type II MOPs associated to $\vec n$ in the system of measures  $\mu_1,\mu_2,\mu_3$ exist if and only if the matrix
    \begin{equation*}
        M_{(1,2,2)} = \left(\begin{array}{c|cc|cc}
            m_{00}^{(1)} & m_{00}^{(2)} & m_{10}^{(2)} & m_{00}^{(3)} & m_{10}^{(3)} \\
            m_{10}^{(1)} & m_{10}^{(2)} & m_{20}^{(2)} & m_{10}^{(3)} & m_{20}^{(3)} \\
            m_{01}^{(1)} & m_{01}^{(2)} & m_{11}^{(2)} & m_{01}^{(3)} & m_{11}^{(3)} \\
            m_{20}^{(1)} & m_{20}^{(2)} & m_{30}^{(2)} & m_{20}^{(3)} & m_{30}^{(3)} \\
            m_{11}^{(1)} & m_{11}^{(2)} & m_{21}^{(2)} & m_{11}^{(3)} & m_{21}^{(3)}
        \end{array}\right)
    \end{equation*}	
    is regular. In that case, $\vec n=(1,2,2)$ would be a normal multi-index.

    Observe that, regarding the shape of the matrix $M_{\vec n}$, defined in \eqref{eq:matrix-Mn} and \eqref{eq:matrix-Mn-blocks}, $M_{(1,2,2)}$ is composed of $r=3$ blocks: the first one of size $5\times 1 = |\vec n|\times n_1$, whereas the remaining blocks size is $5\times 2 = |\vec n|\times n_2 = |\vec n|\times n_3$.
\end{ex}

If we use the polynomial vector notation, it is clear that, for a set of $d+1$ multi-indices $\{\vec n_0,\dots,\vec n_d\}$ satisfying conditions \eqref{eq:condition-multi-indices-1} and \eqref{eq:condition-multi-indices-2}, the existence of bivariate Type II MOPV $\mathbb{P}_{\{\vec n_k\}_0^d}^{(d)}$ is equivalent to the existence of bivariate Type I MOPVs $\mathbb{A}_{\{\vec n_k\}_0^d,j}$ and both are equivalent to the regularity of all matrices $M_{\vec n_0},\dots,M_{\vec n_d}$.

Another basic result in multiple orthogonality, which is a consequence of the definition itself, is a biorthogonality relation satisfied by Type II and Type I MOP (see  \linebreak \citep[Theorem 23.1.6]{Ismail}). We have proven a generalization of this theorem.

\begin{theorem}
    \label{th:biorthogonality-2-app}
    Let $\mu_1,\dots,\mu_r$ be a system of bivariate absolutely continuous measures and consider two normal multi-indices $\vec n,\vec m\in\NN^r$. Then, the bivariate Type II MOP $P_{\vec n}$ and the bivariate Type I function $Q_{\vec m}$ satisfy the following biorthogonality relation:
    \begin{equation}
        \label{eq:biorthogonality-2-app}
        \prodesc{P_{\vec n}}{Q_{\vec m}}_{\mu} = \left\{\begin{array}{ccl}
            0 & \text{ if } & \vec m \leq \vec n \ \text{ (componentwise)}\\
            0 & \text{ if } & |\vec n| \leq |\vec m|-2 \\
            1 & \text{ if } & |\vec n| = |\vec m|-1 
        \end{array}\right.
    \end{equation}
\end{theorem}
\begin{proof}
    Applying the definition of Type I function $Q_{\vec m}$ in \eqref{eq:typeI-function-2-app}, $$\prodesc{P_{\vec n}}{Q_{\vec m}}_{\mu} =\displaystyle\sum_{j=1}^r \prodesc{P_{\vec n}}{A_{\vec m,j}}_j.$$
    \begin{itemize}
        \item If $\vec m \leq \vec n$ (componentwise), then $m_j\leq n_j$ for every $j$. As we know from Definition \ref{def:typeII}, $\pi(\mdeg(A_{\vec m,j}))\leq m_j-1 \leq n_j -1$, due to Type II orthogonality $\prodesc{P_{\vec n}}{A_{\vec m,j}}_j=0$ for every $j$.
        \item If $|\vec n|\leq|\vec m|-2$, then $|\vec n|=\pi(\mdeg(P_{\vec n}))\leq|\vec m|-2$, so that $\displaystyle\sum_{j=1}^r\prodesc{P_{\vec n}}{A_{\vec m,j}}_j = 0$ because of Type I orthogonality.
        \item If $|\vec n|=|\vec m|-1$, then  $|\vec n|=\pi(\mdeg(P_{\vec n}))=|\vec m|-1$ and, again due to Type I orthogonality and because $P_{\vec n}$ is a monic polynomial, $\displaystyle\sum_{j=1}^r\prodesc{P_{\vec n}}{A_{\vec m,j}}_j = 1$.
    \end{itemize}
\end{proof}

If we want to express this orthogonality relation with bivariate MOPV, recall that the result of applying an inner product to two vectors of functions is a matrix whose entries are the scalar inner product applied componentwise:
$$
\prodesc{(f_1,\dots,f_n)}{(g_1,\dots,g_m)^t}_\mu= \left(\prodesc{f_i}{g_j}_\mu\right)_{i,j=1}^{n,m}.
$$

So, the relation is:

\begin{corollary}
    \label{prop:biorthogonality-matrix}
    Consider $\vec n_0,\dots,\vec n_d$ and $\vec m_0,\dots,\vec m_h$, two sets of $d$ and $h$ normal multi-indices satisfying conditions \eqref{eq:condition-multi-indices-1} and \eqref{eq:condition-multi-indices-2}. Then, the bivariate Type II MOPV $\mathbb{P}_{\{\vec n_k\}_0^d}^{(d)}$ defined in Definition \ref{def:typeII-vector} and the bivariate Type I Function $\mathbb{Q}_{\{\vec m_i\}_0^h}$ defined in  \eqref{eq:TypeIVector} satisfy the following relation:
    \begin{equation}
        \label{eq:biorthogonality-matrix}
      \prodesc{\mathbb{P}_{\{\vec n_k\}_0^d}^{(d)}}{\mathbb{Q}_{\{\vec m_i\}_0^h}^t}_\mu = \left\lbrace\begin{array}{ccc}
      0_{(d+1)\times(h+1)}  & \text{ if } & \vec m_h \leq \vec n_0 \\
      \left(\begin{array}{c|c}
      0_{d\times 1} & I_{d} \\ \hline
      0 & 0_{1\times d}
    \end{array}\right) & \text{ if } & \vec m_k = \vec n_k, \ (k=0,\dots,d) \\
      \left(\begin{array}{c|c}
      0_{d\times 1} & 0_{d\times(d+1)} \\ \hline
      1 & 0_{1\times (d+1)}
      \end{array}\right) & \text{ if } & h = d+1 \\
      0_{(d+1)\times(h+1)}  & \text{ if } & h \geq d + 2 \\
      \end{array}\right.
    \end{equation}
\end{corollary}
\begin{proof}
    Just observe
    \begin{equation*}
       \prodesc{\mathbb{P}_{\{\vec n_k\}_0^d}^{(d)}}{\mathbb{Q}_{\{\vec m_i\}_0^h}^t}_\mu = \begin{pmatrix}
\prodesc{P_{\vec n_0}}{Q_{\vec m_0}} & \cdots & \prodesc{P_{\vec n_0}}{Q_{\vec m_h}} \\
\vdots & \ddots & \vdots \\
\prodesc{P_{\vec n_d}}{Q_{\vec m_0}} & \cdots & \prodesc{P_{\vec n_d}}{Q_{\vec m_h}} 
\end{pmatrix} 
    \end{equation*}
    and apply \eqref{eq:biorthogonality-2-app} componentwise.

\end{proof}

Equation \eqref{eq:biorthogonality-matrix} may appear quite asymmetric. This is a consequence of the asymmetry of equation \eqref{eq:biorthogonality-2-app}, where the result is $1$ if $|\vec n|=|\vec m|-1$. Even the original univariate biorthogonality relation in \citep{Ismail} shares this lack of symmetry.

The next relations we are going to generalize are the `Nearest Neighbours Recurrence Relations', which are in turn generalizations of the three-term relation by using MOPs (see \citep[Theorem 23.1.7, Theorem 23.1.9]{Ismail}).

\begin{figure}[h]
    \centering\includegraphics[width=5cm]{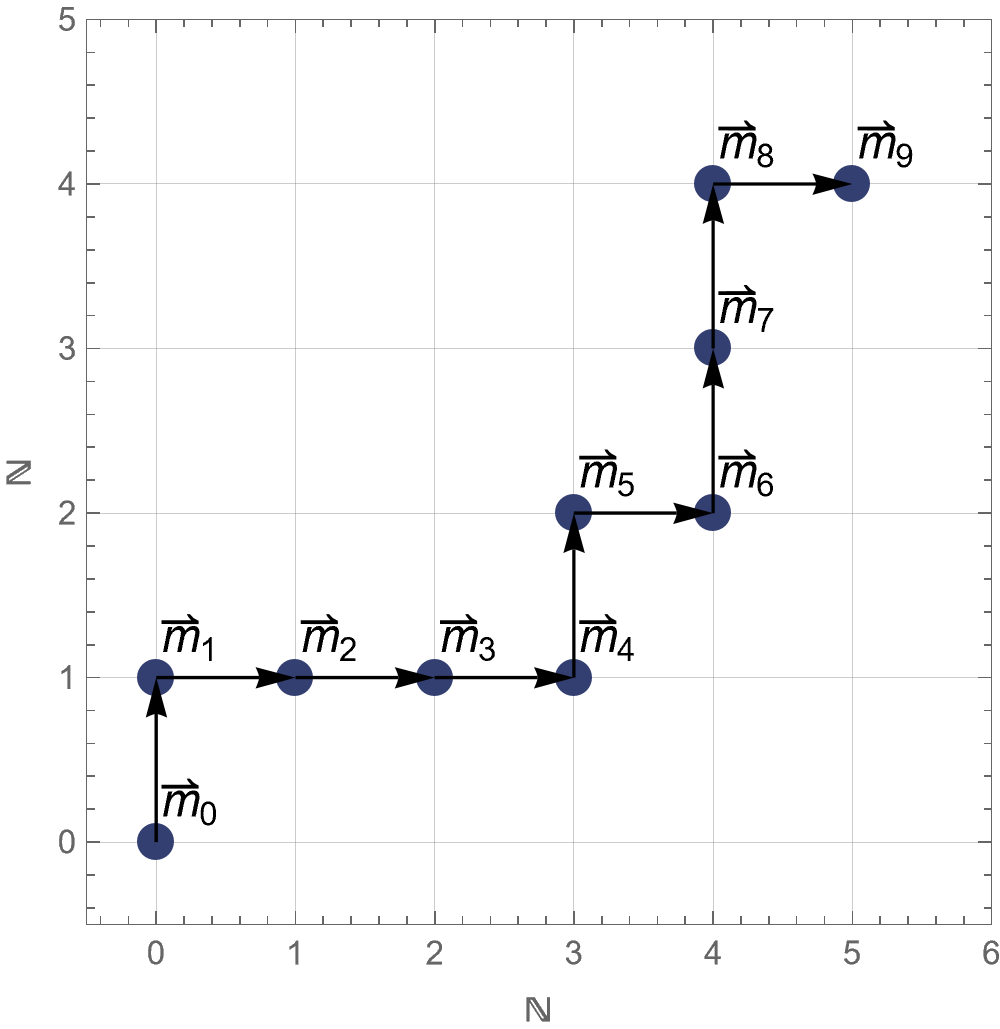}
    \caption{A path from $(0,0)$ to $(5,4)$ for $r=2$}
    \label{fig:path}
\end{figure}

Since we are working with multi-indices and bivariate polynomials, there are several ways to decrease or increase the degree or the multidegree of a bivariate MOP. For example, given two multi-indices $\vec n, \vec w$ such that $\vec n\leq \vec w$ componentwise, we know that $$|\vec n|=\pi(\mdeg(P_{\vec n}))\leq |\vec w|=\pi(\mdeg(P_{\vec w})),$$ and as a consequence, $\deg(P_{\vec n})\leq\deg(P_{\vec w})$. Along the remainder of this section, we will consider paths $\{\vec m_k: k=0,\dots,|\vec n|,\dots,|\vec w|\}$ from $\vec 0=(0,\dots,0)$ to $\vec w$ via $\vec n$, with $\vec m_0=\vec 0$, $\vec m_{|\vec n|}=\vec n$, $\vec m_{|\vec w|}=\vec w$ and where in each step the multi-index $\vec m_k$ is increased by one at exactly one component, \textit{i.e.}, the path is composed of neighbour multi-indices (they satisfy condition \eqref{eq:condition-multi-indices-2} $k=0,\dots,|\vec w|-1$). For these paths, we have $|\vec m_k| = k$ and $\vec m_k\leq \vec m_{k+1}$. See Figure \ref{fig:path} for an example of these kind of path.

Whereas the bivariate biorthogonality relation is almost the same as the univariate one in \citep[Theorem 23.1.6]{Ismail}, we will show that the extended Nearest Neighbour Relations (NNR) are different; however, this issue can be solved by using the vector notation. There are two different types of NNR, one employs bivariate Type II MOP while the other one uses bivariate Type I MOP (Type I functions, actually). 

Firstly, we will focus on the NNR for bivariate Type II MOP. Given $\vec n\in\NN^r$, the original univariate Nearest Neighbours Relation \citep[Theorem 23.1.7]{Ismail} expresses the polynomial $xP_{\vec n}(x)$ as a linear combination of $P_{\vec n}$ itself, a neighbour of higher degree and $r$ neighbours of lower degree. In this case, it is clear that $\deg(x P_{\vec n})=\deg(P_{\vec n})+1=|\vec n|+1$, and that is the reason why only one neighbour of higher degree is needed. However, this equality does not hold in the bivariate case, since $\pi(\mdeg(xP_{\vec n}))\neq \pi(\mdeg(P_{\vec n})) +1$, and something similar happens when we multiply by `$y$'. For this reason, we present the following lemma.

\begin{lemma}
    \label{lemma:pi-xyP}
    Let 
    $$P(x,y)= \underbrace{c_{\pi(l,m)} x^l y^m}_{\text{leading term}}+\cdots + c_1 x + c_0, \qquad c_{\pi(l,m)}\neq 0$$ 
    so that $(l,m)=\mdeg(P(x,y))$. Then, 
    \begin{equation*}
        \begin{array}{c}
            \pi(l+1,m)=\pi(\mdeg(x P(x,y)))=\pi(\mdeg(P(x,y)))+\deg(P(x,y))+1, \\
            \pi(l,m+1)=\pi(\mdeg(y P(x,y)))=\pi(\mdeg(P(x,y)))+\deg(P(x,y))+2,		
        \end{array}
    \end{equation*}
    where $\deg(P(x,y)) = l+m$ is the total degree of $P(x,y)$. 
\end{lemma}
\begin{proof}
    Remember the usual basis of the vector space of bivariate polynomials ordered by the reverse lexicographic graded order \eqref{eq:usual-basis-ordered}:
    $$
    \{1,\dots, x^{l+m},\dots, x^l y^m, \underbrace{x^{l-1}y^{m+1}, \dots, x y^{l+m-1},y^{l+m}}_{l\text{ monomials}},\underbrace{x^{l+m+1},\dots,x^{l+2}y^{m-1},x^{l+1}, y^m}_{m\text{ monomials}},x^l y^{m+1},\dots\}.
    $$
    An alternative interpretation of the function $\pi(l,m)$ is the position of the monomial $x^l y^m$ on this list (starting the count in $0$). In this way, there are $l+m$ monomials in between $x^l y^m$ and $x^{l+1} y^m$. Thus, $\pi(l+1,m)=\pi(l,m)+l+m+1 = \pi(\mdeg(P)) + \deg(P) + 1$. For the second equality, observe the monomial $x^l y^{m+1}$ is exactly the next one after $x^{l+1}y^m$.
\end{proof}

With this, we can present different results of NNR for Type II bivariate MOP.

\begin{theorem}[Nearest Neighbour Relation (for $xP_{\vec n}$ and $y P_{\vec n}$)]
    \label{th:NNR-xyPn}
    Let us consider $(\mu_1,\dots,\mu_r)$, a perfect system of bidimensional measures, a multi-index $\vec n=(n_1,\dots,n_r)\in\NN^r$, $P_{\vec n}$ its associated bivariate Type II MOP and $d_{\vec n}, k_{\vec n}$ its degree and remainder introduced in Definition \ref{def:dnkn}. Define $\vec v = (n_1-d_{\vec n}-1,\dots,n_r-d_{\vec n}-1)$. Choose $\vec w_x, \vec w_y \in \NN^r$ such that 
        \begin{align*}
            |\vec w_x|&=\frac 1 2 (d_{\vec n}+1)(d_{\vec n}+2) + k_{\vec n} = |\vec n| + d_{\vec n}+1, & |\vec w_y|&=\frac 1 2 (d_{\vec n}+1)(d_{\vec n}+2) + k_{\vec n}+1 = |\vec n| + d_{\vec n}+2
        \end{align*}
    and $\vec w_y\geq \vec w_x\geq \vec v$ componentwise, and let $P_{\vec w_x}, P_{\vec w_y}$ be their associated bivariate Type II MOP. Consider $\{\vec m_i: i=|\vec n|-(d_{\vec n}+1)r,\dots,|\vec n| + d_{\vec n}+2\}$ a path of neighbour multi-indices where $|\vec m_i| = i$ for every $i$, $\vec m_{|\vec n|-(d_{\vec n}+1)r}=\vec v$, $\vec m_{|\vec n|}=\vec n$, $\vec m_{|\vec n| + d_{\vec n}+1}=\vec w_x$ and $\vec m_{|\vec n| + d_{\vec n}+2}=\vec w_y$. For every $i$, let $P_{\vec m_i}$ be the bivariate Type II MOP and  $Q_{\vec m_i}$ the bivariate Type I function \eqref{eq:typeI-function-2-app} associated to the multi-index $\vec m_i$. Then, the following relations hold:
    \begin{equation}
        \label{eq:NNR-xPn}
        x P_{\vec n}(x,y) = P_{\vec w_x}(x,y) + \sum_{i=|\vec n|-(d_{\vec n}+1)r}^{|\vec n| + d_{\vec n}} a_i P_{\vec m_i}(x,y),
    \end{equation}
    \begin{equation}
        \label{eq:NNR-yPn}
        y P_{\vec n}(x,y) = P_{\vec w_y}(x,y) + \sum_{i=|\vec n|-(d_{\vec n}+1)r}^{|\vec n| + d_{\vec n}+1} b_i P_{\vec m_i}(x,y),
    \end{equation}
    where $a_{i-1}=\prodesc{xP_{\vec n}(x,y)}{Q_{\vec m_i}(x,y)}$ and $b_{i-1}=\prodesc{yP_{\vec n}(x,y)}{Q_{\vec m_i}(x,y)}$.
    
\end{theorem}
\begin{proof}
We will prove \eqref{eq:NNR-xPn}. The proof of \eqref{eq:NNR-yPn} is analogous but using the second equation of Lemma \ref{lemma:pi-xyP} instead of the first one when necessary. 

Let us consider $\{\vec m_i: i=0,\dots,|\vec n| + d_{\vec n}+1\}$ a path of neighbour multi-indices where $|\vec m_i| = i$ $(i=0,\dots,|\vec n| + d_{\vec n}+1)$, $\vec m_0 =\vec 0$, $\vec m_{|\vec n|-(d_{\vec n}+1)r}=\vec v$, $\vec m_{|\vec n|}=\vec n$ and $\vec m_{|\vec n| + d_{\vec n}+1}=\vec w$. 

According to Lemma \ref{lemma:pi-xyP}, we know that $\pi(\mdeg(xP_{\vec n}))=|\vec n|+d_{\vec n}+1=|\vec w|$. In addition, as $P_{\vec n}$ is monic, so is $xP_{\vec n}$. On the other hand, as $\pi(\mdeg(P_{\vec m_k}))=k$, then $\{P_{\vec m_0}, P_{\vec m_1} \dots,P_{\vec m_{|\vec n|+d_{\vec n}+1}}\}$ is a basis of the vector space of bivariate polynomials $P(x,y)$ such that \linebreak $\pi(\mdeg(P))\leq |\vec n|+d_{\vec n}+1 = |\vec w|$. With all this information, it is possible to write $xP_{\vec n}$ as a linear combination of $P_{\vec m_0}, P_{\vec m_1} \dots,P_{\vec m_{|\vec n|+d_{\vec n}+1}}$.
\begin{equation}
    xP_{\vec n} =  P_{\vec w_x} + \sum_{k=0}^{|\vec n|+d_{\vec n}}a_k P_{\vec m_k}.
\end{equation}
In order to find the values of coefficients $a_k$, for $k=0,\dots,|\vec n|+d_{\vec n}$, observe
\begin{equation}
    \label{eq:NNR-2-app-1}
    \prodesc{xP_{\vec n}}{Q_{\vec m_l}} = \prodesc{P_{\vec w}}{Q_{\vec m_l}} + \sum_{k=0}^{|\vec n|+d_{\vec n}} a_k \prodesc{P_{\vec m_k}}{Q_{\vec m_l}}.
\end{equation}
To simplify the previous expression, we will employ the biorthogonality relation of Theorem \ref{th:biorthogonality-2-app}. 
\begin{itemize}
    \item Firstly, it is known that $\vec m_l \leq \vec m_k$ if, and only if $l\leq k$, then $\prodesc{P_{\vec m_k}}{Q_{\vec m_l}}=0$ if $k\geq l$.
    \item Secondly, because of $|\vec m_k|=k$,  $|\vec m_k| \leq |\vec m_l| - 2$ if and only if $k\leq l-2$, then $\prodesc{P_{\vec m_k}}{Q_{\vec m_l}}=0$ if $k\leq l-2$.
    \item Lastly, if $k=l-1$, then $\prodesc{P_{\vec m_{l-1}}}{Q_{\vec m_l}}=1$.
\end{itemize}
Applying these conclusions to \eqref{eq:NNR-2-app-1}, we get
\begin{equation}
    \prodesc{xP_{\vec n}}{Q_{\vec m_l}}=a_{l-1}, \ \ \ l=1,\dots,|\vec n|+d_{\vec n}+1
\end{equation}
The next step is to prove that $a_0=\dots=a_{|\vec n|-(d_{\vec n}+1)r-1}=0$, \textit{i.e.}, $a_{k-1} = \prodesc{xP_{\vec n}}{Q_{\vec m_k}} = 0$ for $k=1,\dots,|\vec n|-(d_{\vec n}+1)r$. Applying \eqref{eq:typeI-function-2-app}:
$$
\prodesc{xP_{\vec n}}{Q_{\vec m_k}} = \sum_{j=1}^r \prodesc{xP_{\vec n}}{A_{\vec m_k,j}}_j,
$$
and because of the integral definition of these inner products 
$$ \displaystyle\sum_{j=1}^r \prodesc{xP_{\vec n}}{A_{\vec m_k,j}}_j=\displaystyle\sum_{j=0}^r \prodesc{P_{\vec n}}{xA_{\vec m_k,j}}_j.$$
Fixing $j\in\{1,\dots,r\}$, and due to Type II orthogonality, we know that $$\prodesc{P_{\vec n}}{xA_{\vec m_k,j}}_j=0\qquad \text{ if } \pi(\mdeg(xA_{\vec m_k,j}))\leq n_j-1.$$
Now, remember $k\leq |\vec n|-(d_{\vec n}+1)r$ and, in this case, $\vec m_k \leq (n_1-d_{\vec n}-1,\dots,n_r-d_{\vec n}-1)=\vec v$, so that necessarily, $\pi(\mdeg(A_{\vec m_k,j}))\leq n_j-d_{\vec n}-2$ and $\deg(A_{\vec m_k,j}) \leq d_{\vec n}-1$. With these inequalities and using Lemma \ref{lemma:pi-xyP} we have
$$
\pi(\mdeg(x A_{\vec m_k,j} )) = \pi(\mdeg(A_{\vec m_k,j})) + \deg(A_{\vec m_k,j}) + 1 \leq n_j - d_{\vec n}-2 + d_{\vec n} - 1 + 1 \leq n_j-1.
$$
Then, $\prodesc{P_{\vec n}}{xA_{\vec m_k,j}}_j=0$ for $j=1,\dots,r$ and $a_{k-1} = \prodesc{xP_{\vec n}}{Q_{\vec m_k}} = 0$ for  $k=1,\dots,|\vec n|-(d_{\vec n}+1)r$.

\end{proof}

Due to the little difference in the behaviour of the variables `$x$' and `$y$' shown in Lemma~\ref{lemma:pi-xyP}, we also have some differences in the items of this last theorem. On one hand, in \eqref{eq:NNR-yPn} there are $d_{\vec n}+2$ multi-indices of higher norm, one more than in \eqref{eq:NNR-xPn}. On the other hand, in both items we employ $(d_{\vec n}+1)r$ multi-indices of lower norm. To clarify this result, we present an example.

\begin{ex}
    Assume $r=2$, $\vec n =(6,8)$, then $|\vec n|=14=\pi(0,4)$, $d_{\vec n}=4$, $k_{\vec n}=4$, \linebreak$\vec v = (n_1-d_{\vec n}-1,n_2-d_{\vec n}-1)=(1,3)$, and we need $\vec w_y\geq\vec w_x\geq \vec n$ s.t. $|\vec w_x| = 5\cdot 6 /2 + 4 = \pi(1,4) = 19$ and $|\vec w_y| = 5\cdot 6 /2 + 4 +1=\pi(0,5)= 20$, let us take, for example $\vec w_x = (9,10)$ and $w_y = (9,11)$. Now, we choose a path from $\vec 0$ to $\vec w_y$ via $\vec v$, $\vec n$ and $\vec w_x$.
    $$
    \vec m_0 =(0,0)\rightarrow \vec m_1 =(0,1)\rightarrow \vec m_2 =(0,2)\rightarrow \vec m_3 =(1,2)\rightarrow \underline{\vec v = \vec m_4 =(1,3)}\rightarrow 
    $$ $$
    \rightarrow\vec m_5 = (2,3)\rightarrow \vec m_6 =(2,4)\rightarrow \vec m_7 =(2,5)\rightarrow \vec m_8 = (2,6)\rightarrow \vec m_9 =(3,6)\rightarrow
    $$ $$
    \rightarrow\vec m_{10} = (3,7)\rightarrow \vec m_{11} =(4,7)\rightarrow \vec m_{12} =(5,7)\rightarrow \vec m_{13} = (6,7)\rightarrow \underline{\vec n = \vec m_{14} =(6,8)}\rightarrow
    $$ $$\rightarrow \vec m_{15} = (7,8)\rightarrow \vec m_{16} = (7,9)\rightarrow \vec m_{17} = (8,9)\rightarrow \vec m_{18} = (9,9)\rightarrow$$ $$\rightarrow \underline{\vec w_x = \vec m_{19} = (9,10)}\rightarrow \underline{\vec w_y = \vec m_{20} = (9,11)}.$$
    Then
    \begin{align}
        \label{eq:example-NNR}
    xP_{(6,8)} &= P_{(9,10)} + \displaystyle\sum_{i=4}^{18} a_i P_{\vec m_i}, &
    yP_{(6,8)} &= P_{(9,11)} + \displaystyle\sum_{i=4}^{19} b_i P_{\vec m_i}. 
    \end{align}

    Observe Figure \ref{fig:example-NNR}(a), where we represent the multi-indices employed in the linear combinations \eqref{eq:example-NNR}. As shown, $5=d_{\vec n}+1$ multi-indices whose norms are higher than $14=|\vec n|$ and $10= (d_{\vec n}+1)r$ multi-indices whose norms are lower than $|\vec n|$ are used. $P_{\vec m_0}, P_{\vec m_1}, P_{\vec m_2}$ and $P_{\vec m_3}$ do not take place in \eqref{eq:example-NNR}.

    Additionally, the relation still holds if we choose a different path from $\vec v$ to $\vec n$, and even if we choose another $\vec w_x$ or $\vec w_y$. For example, considering the following path (see Figure \ref{fig:example-NNR}(b)):

    $$
    \vec m_0 =(0,0)\rightarrow \vec m_1 =(1,0)\rightarrow \vec m_2 =(1,1)\rightarrow \vec m_3 =(1,2)\rightarrow \underline{\vec v = \vec m_4 =(1,3)}\rightarrow 
    $$ $$
    \rightarrow\vec m_5 = (1,4)\rightarrow \vec m_6 =(1,5)\rightarrow \vec m_7 =(1,6)\rightarrow \vec m_8 = (2,6)\rightarrow \vec m_9 =(3,6)\rightarrow
    $$ $$
    \rightarrow\vec m_{10} = (3,7)\rightarrow \vec m_{11} =(3,8)\rightarrow \vec m_{12} =(4,8)\rightarrow \vec m_{13} = (5,8)\rightarrow \underline{\vec n = \vec m_{14} =(6,8)}\rightarrow
    $$ $$\rightarrow \vec m_{15} = (7,8)\rightarrow \vec m_{16} = (8,8)\rightarrow \vec m_{17} = (9,8)\rightarrow \vec m_{18} = (10,8)\rightarrow $$ $$\rightarrow \underline{\vec w_x = \vec m_{19} = (11,8)}\rightarrow \underline{\vec w_y = \vec m_{20} = (12,8)}.$$
    Then
    \begin{align}
        \label{eq:example-NNR-2}
    xP_{(6,8)} &= P_{(11,8)} + \displaystyle\sum_{i=4}^{18} a_i P_{\vec m_i}, & yP_{(6,8)} &= P_{(12,8)} + \displaystyle\sum_{i=4}^{19} b_i P_{\vec m_i}.
    \end{align}

    \begin{figure}[h]
        \centering\begin{tabular}{cc}
            \includegraphics[width=5cm]{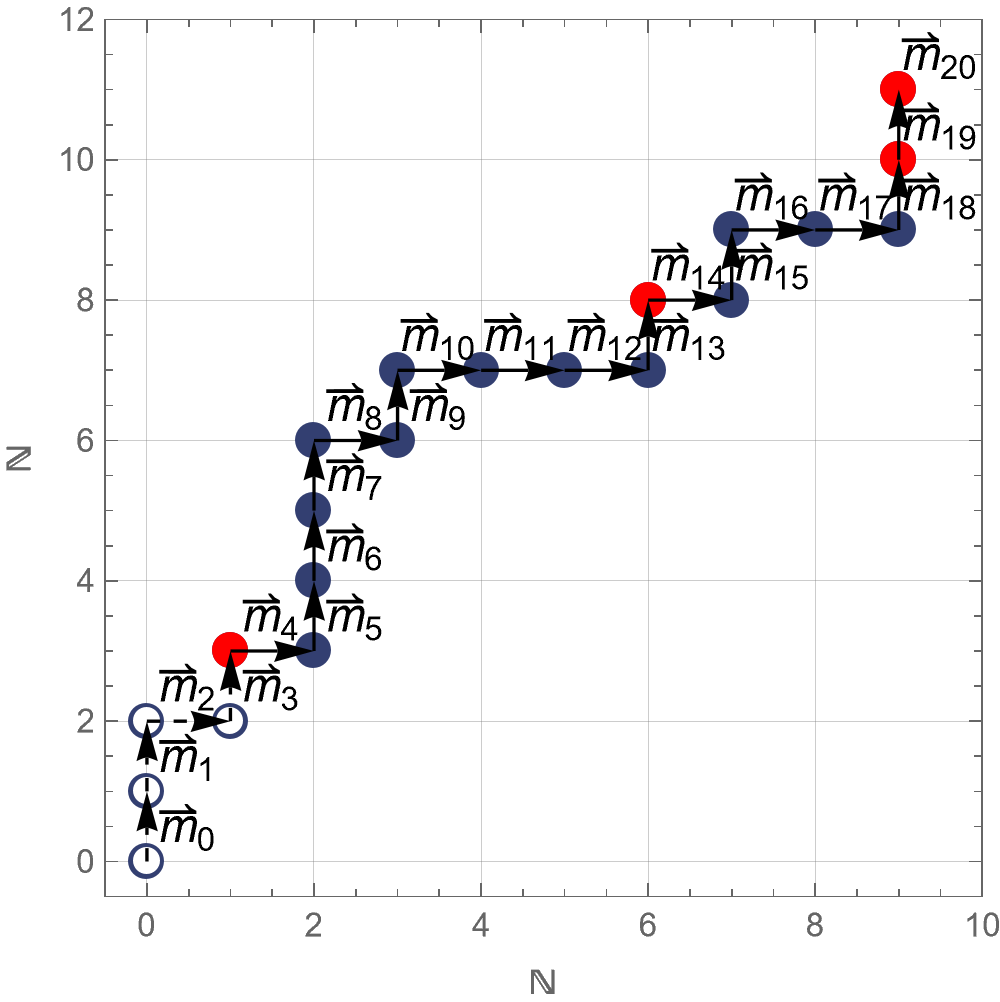} & \includegraphics[width=5cm]{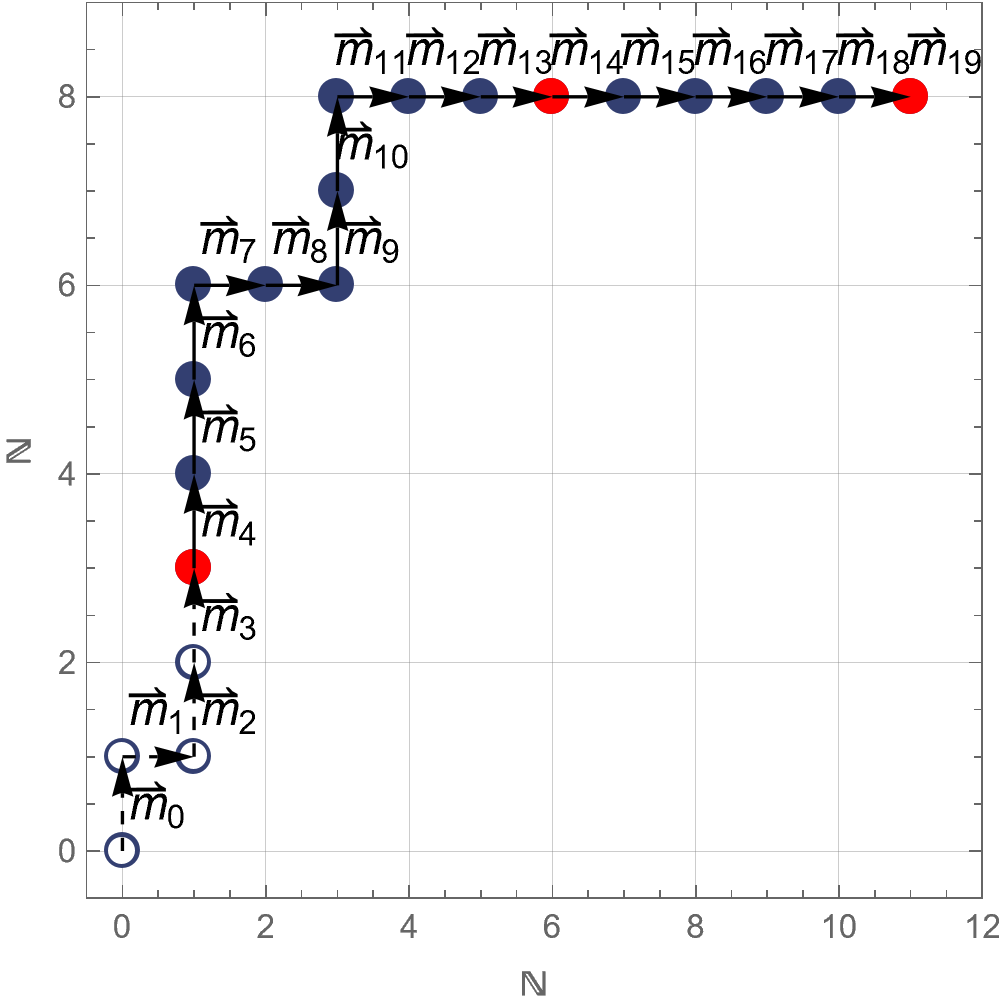} \\
            (a) & (b)
        \end{tabular}
        \caption{Example of NNR for bivariate Type II MOP}
        \label{fig:example-NNR}
    \end{figure}
\end{ex}

When using the NNR, the polynomial vector notation is really useful for clarifying the relation and its notation. Remember a proper bivariate Type II MOPV $\mathbb{P}_{\{\vec n_k\}_0^d}^{(d)}$ is composed of $d+1$ polynomials whose leading terms are, respectively, $x^d,x^{d-1}y,\dots,y^d$, so that $\pi(\mdeg(P_{\vec n_k})) = \frac 1 2 d(d+1) + k$. So, if we bring all the polynomials of equal degree together in polynomial vectors, it is possible to obtain a vector version of the NNR.

\begin{corollary}[Nearest Neighbour Relation (for $x \mathbb{P}_{\{\vec n_k\}_0^d}^{(d)}$ and $y \mathbb{P}_{\{\vec n_k\}_0^d}^{(d)}$)] Let $\mu_1,\dots,\mu_r$ be a perfect system of $r$ bidimensional measures. Consider $d\in\NN$ and $d+1$ multi-indices $\vec n_0,\dots,\vec n_d$ satisfying \eqref{eq:condition-multi-indices-1} and \eqref{eq:condition-multi-indices-2} and their associated bivariate Type II MOPV $\mathbb{P}_{\{\vec n_k\}_0^d}^{(d)}$. Define 
    $$K:= \max\left\{h\in\{0,\dots,d-1\}: \frac 1 2 h(h+1) < |\vec n_0| -(d+1)r\right\}$$
and let $\{\vec m_{h,i},\ \ h=K,\dots,d-1, \ i=0,\dots,h\}$ be a path of neighbour multi-indices with $|\vec m_{h,i}| = \frac 1 2 h(h+1) + i$, where the last multi-index $\vec m_{d-1,d-1}$ is a neighbour of $\vec n_0$ and the multi-index $(\vec n_{0,1}-d-1,\dots,\vec n_{0,r}-d-1)$ belongs to the path. This means, there exists $0\leq l\leq K$ with
$$
\vec m_{K,l}=(\vec n_{0,1}-d-1,\dots,\vec n_{0,r}-d-1).
$$
Let $\mathbb{P}_{\{\vec m_{h,i}\}_0^h}^{(h)}$ be the bivariate Type II MOPV associated to the set of multi-indices $\{\vec m_{h,0},\dots,\vec m_{h,h}\}$, \ \ $h=K,\dots,d-1$. Choose $\vec v_0,\dots,\vec  v_{d+1}$ satisfying \eqref{eq:condition-multi-indices-1} and \eqref{eq:condition-multi-indices-2}, and $\vec n_d$ is a neighbour of $\vec v_0$. Let $\mathbb{P}_{\{\vec v_k\}_0^{d+1}}^{(d+1)}$ be the bivariate Type II MOPV associated.
    Then
    \begin{equation*}
      x \mathbb{P}_{\{\vec n_k\}_0^d}^{(d)} = L_{d+1,x} \mathbb{P}_{\{\vec v_k\}_0^{d+1}}^{(d+1)} + A_d \mathbb{P}_{\{\vec n_k\}_0^d}^{(d)} + \sum_{h=K-1}^{d-1} A_h \mathbb{P}_{\{\vec m_{h,i}\}_0^h}^{(h)}
    \end{equation*}
    and 
    \begin{equation*}
        y \mathbb{P}_{\{\vec n_k\}_0^d}^{(d)} = L_{d+1,y} \mathbb{P}_{\{\vec v_k\}_0^{d+1}}^{(d+1)} + B_d \mathbb{P}_{\{\vec n_k\}_0^d}^{(d)} + \sum_{h=K-1}^{d-1} B_h \mathbb{P}_{\{\vec m_{h,i}\}_0^h}^{(h)}
    \end{equation*}
    where $A_h$ and $B_h$ are matrices of dimension $(d+1)\times(h+1)$, $L_{d+1,x}=(I_{d+1} \ | \ 0_{(d+1)\times 1})$ and $L_{d+1,y}=(0_{(d+1)\times 1} \ | \ I_{d+1})$.
  \end{corollary}
  \begin{proof}
    In this proof, only the relation for $x \mathbb{P}_{\{\vec n_k\}_0^d}^{(d)}$ will be proven. The arguments for $y \mathbb{P}_{\{\vec n_k\}_0^d}^{(d)}$ are analogous but applying the second equation of Lemma \ref{lemma:pi-xyP}.

    Similarly to the proof of Theorem \ref{th:NNR-xyPn}, let us consider a path of neighbour multi-indices 
    $$\{\vec m_{h,i}: h=0,\dots,d-1, \ \ i=0,\dots,h \} = \{\vec m_{0,0},\vec m_{1,0},\vec m_{1,1},\vec m_{2,0},\vec m_{2,1},\vec m_{2,2},\dots,\vec m_{d-1,0},\dots,\vec m_{d-1,d-1}\}$$
    so that $|\vec m_{h,i}| = \frac 1 2 h(h+1) + i$, $\vec m_{K,l}=(\vec n_{0,1}-d-1,\dots,\vec n_{0,r}-d-1)$ and $\vec m_{d-1,d-1}\leq \vec n_0$.

    For $h=0,\dots,d-1$, we know the polynomial vectors $\mathbb{P}_{\{\vec m_{h,i}\}_0^h}^{(h)}$ are composed of $h+1$ polynomials of degree $h$ and their leading terms are $x^h,x^{h-1}y,\dots,y^h$. Furthermore, the polynomial vector $\mathbb{P}_{\{\vec n_k\}_0^{d}}^{(d)}$ (resp. $\mathbb{P}_{\{\vec v_k\}_0^{d+1}}^{(d+1)}$) is composed of $d+1$ (resp. $d+2$) polynomials of degree $d$ (resp. $d+1$) and their leading terms are $x^d,x^{d-1}y,\dots,y^d$ (resp. $x^{d+1},x^{d}y,\dots,y^{d+1}$). Then, the polynomial vectors
    $$
    \left\{ \mathbb{P}_{\{\vec m_{0,i}\}_0^0}^{(0)},\mathbb{P}_{\{\vec m_{1,i}\}_0^1}^{(1)},\dots,\mathbb{P}_{\{\vec m_{d-1,i}\}_0^{d-1}}^{(d-1)}, \mathbb{P}_{\{\vec n_k\}_0^{d}}^{(d)}, \mathbb{P}_{\{\vec v_k\}_0^{d+1}}^{(d+1)}\right\}
    $$
    shape a basis of the space of polynomial vectors of degree $\leq d+1$ (where the coefficients of the linear combinations are matrices of proper dimension). It is clear that $x \mathbb{P}_{\{\vec n_k\}_0^d}^{(d)}$ is a polynomial vector composed of $d+1$ polynomials of degree $d+1$. Then, it is possible to express
    \begin{equation}
        \label{eq:NNRV-1}
        x \mathbb{P}_{\{\vec n_k\}_0^d}^{(d)} = A_{d+1}\mathbb{P}_{\{\vec v_k\}_0^{d+1}}^{(d+1)} + A_d \mathbb{P}_{\{\vec n_k\}_0^{d}}^{(d)} + \sum_{h=0}^{d-1} A_h \mathbb{P}_{\{\vec m_{h,i}\}_0^h}^{(h)},
    \end{equation}
    where $A_h$ are matrices of dimension $(d+1)\times(h+1)$ for $h=0,\dots,d+1$.
    
    Observe that the leading terms of $x \mathbb{P}_{\{\vec n_k\}_0^d}^{(d)}$ (resp. $y \mathbb{P}_{\{\vec n_k\}_0^d}^{(d)}$) are the same as the leading terms of the first (resp. last) $d+1$ polynomials of $\mathbb{P}_{\{\vec v_k\}_0^{d+1}}^{(d+1)}$. Then, necessarily $A_{d+1}=L_{d+1,x}$ (resp. $B_{d+1}=L_{d+1,y}$).

    Now, we want to show that $A_0=A_1=\cdots=A_{K-2}=0$. First, if we multiply \eqref{eq:NNRV-1} by the Type I Function $\mathbb{Q}_{\{\vec m_{t,s}\}_0^t}$ we get
    \begin{equation}
        \prodesc{x \mathbb{P}_{\{\vec n_k\}_0^d}^{(d)}}{\mathbb{Q}_{\{\vec m_{t,s}\}_0^t}} = L_{d,x}\prodesc{\mathbb{P}_{\{\vec v_k\}_0^{d+1}}^{(d+1)}}{\mathbb{Q}_{\{\vec m_{t,s}\}_0^t}} + A_d \prodesc{\mathbb{P}_{\{\vec n_k\}_0^{d}}^{(d)}}{\mathbb{Q}_{\{\vec m_{t,s}\}_0^t}} + \sum_{h=0}^{d-1} A_h \prodesc{\mathbb{P}_{\{\vec m_{h,i}\}_0^h}^{(h)}}{\mathbb{Q}_{\{\vec m_{t,s}\}_0^t}},
    \end{equation}
    Applying the biorthogonality relation for bivariate MOPV \eqref{eq:biorthogonality-matrix}, we get that, for a fixed $t$:
    $$
    \prodesc{x \mathbb{P}_{\{\vec n_k\}_0^d}^{(d)}}{\mathbb{Q}_{\{\vec m_{t,s}\}_0^t}} = A_{t} \left(\begin{array}{c|c}
        0_{t\times 1} & I_{t} \\ \hline
        0 & 0_{1\times t} 
    \end{array}\right) + A_{t-1} \left(\begin{array}{c|c}
            0_{(t-1) \times 1} & 0_{(t-1)\times t} \\ \hline
            1 & 0_{1\times t}
            \end{array}\right).
    $$
    Then, observe the first summand is a matrix whose first column is $0$ and the other $t$ columns are the first $t$ columns of $A_t$. On the other hand, the second summand is a matrix where the first column is the last column of $A_{t-1}$ and the remaining columns are $0$. So, we have that $\prodesc{x \mathbb{P}_{\{\vec n_k\}_0^d}^{(d)}}{\mathbb{Q}_{\{\vec m_{t,s}\}_0^t}} $ is equal to a matrix whose first column is the last column of $A_{t-1}$ and the other $t$ columns are the first $t$ columns of $A_t$. This means:
    \begin{equation}
        \label{eq:NNRV-3}
        \prodesc{x\mathbb{P}_{\{\vec n_k\}_0^d}^{(d)}}{(\mathbb{Q}_{\{\vec m_{1,s}\}_0^1}|\mathbb{Q}_{\{\vec m_{2,s}\}_0^2}| \ \cdots \ |\mathbb{Q}_{\{\vec m_{d+1,s}\}_0^{d+1}})}=(A_0|A_1| \ \cdots \ |A_d|I_{d+1}).
    \end{equation}
    So, we have to ensure that, for $1\leq h \leq K-1$, $\prodesc{x \mathbb{P}_{\{\vec n_k\}_0^d}^{(d)}}{\mathbb{Q}_{\{\vec m_{h,l}\}_0^h}} = 0$. We have 

    \begin{equation}
        \prodesc{x \mathbb{P}_{\{\vec n_k\}_0^d}^{(d)}}{\mathbb{Q}_{\{\vec m_{h,l}\}_0^h}} = \sum_{j=1}^r \prodesc{ \mathbb{P}_{\{\vec n_k\}_0^d}^{(d)}}{x \mathbb{A}_{\{\vec m_{h,l}\}_0^h,j}}_j.
    \end{equation}
    $\prodesc{ \mathbb{P}_{\{\vec n_k\}_0^d}^{(d)}}{x \mathbb{A}_{\{\vec m_{h,l}\}_0^h,j}}_j = 0 $ if 
    $$
    \pi(\mdeg(xA_{(\vec m_{h,l}),j}))\leq \vec n_{k,j}-1 \ \ \ \forall \ \ l=0,\dots,h, \ \ k=0,\dots, d, \ \ j=1,\dots,r . 
    $$ Remember all these multi-indices are neighbour, so that the previous condition is equivalent to   
    \begin{equation}
        \label{eq:NNRV-2}
        \pi(\mdeg(xA_{(\vec m_{h,h}),j}))\leq \vec n_{0,j}-1.
    \end{equation}
    But, as $h \leq K-1$, necessarily $\vec m_{h,h} \leq (n_{0,1}-d-1,\dots,n_{0,r}-d-1)$. Applying the same argument as in the proof of Theorem \ref{th:NNR-xyPn}, we get \eqref{eq:NNRV-2}, so that $\prodesc{x \mathbb{P}_{\{\vec n_k\}_0^d}^{(d)}}{\mathbb{Q}_{\{\vec m_{h,l}\}_0^h}}=0$ for $h=1,\dots,K-1$. Observing \eqref{eq:NNRV-3}, this means $A_0=A_1=\dots=A_{K-2}=0$, and $A_{K-1}$ is a matrix whose first $d+1$ columns are $0$.

  \end{proof}

  \begin{remark}
    This proof could also be deduced by applying Theorem \ref{th:NNR-xyPn} to every single polynomial of the vector $x\mathbb{P}_{\{\vec n_k\}_0^d}^{(d)}$ and reuniting the coefficients in the matrices $A_h$,\linebreak $h=K-1,\dots,d+1$. However, we present this one in order to apply the definitions and results with polynomial vectors.
  \end{remark}

  Observe this relation looks really similar to the univariate \citep[Theorem 23.1.7]{Ismail}.

   Previously in this section we mentioned that in the univariate case there exists a NNR for Type I MOP, summarized employing the Type I function, see \citep[Theorem 23.1.9]{Ismail}. We have also proven the extension of this result for bivariate MOP.

    \begin{theorem}[Nearest Neighbour Relation (for $x Q_{\vec n}$ and $y Q_{\vec n}$)]
        \label{th:NNR-xyQn}
        Let us \linebreak consider $(\mu_1,\dots,\mu_r)$, a perfect system of bidimensional and absolutely continuous measures, a multi-index $\vec n=(n_1,\dots,n_r)\in\NN^r$, its associated degree and remainder $d_{\vec n}, k_{\vec n}$ introduced in Definition \ref{def:dnkn} and its associated bivariate Type I function $Q_{\vec n}$ defined in \eqref{eq:typeI-function-2-app}.
        \begin{enumerate}
            \item Consider $\{\vec m_k: k=|\vec n|-d_{\vec n},\dots,|\vec n| + (d_{\vec n}+1)r\}$ a path of neighbour multi-indices where $|\vec m_k| = k$, $\vec m_{|\vec n|}=\vec n$ and $\vec m_{|\vec n| + (d_{\vec n}+1)r}=(n_1+d_{\vec n}+1,\dots,n_r+d_{\vec n}+1)$. For every k, let $Q_{\vec m_k}$ be the bivariate Type I Function and $P_{\vec m_k}$ the bivariate Type II MOP. Then:
            $$
            x Q_{\vec n}(x,y) = Q_{\vec m_{|\vec n|-d_{\vec n}}}(x,y) + \sum_{k=|\vec n|-d_{\vec n}+1}^{|\vec n|+(d_{\vec n}+1)r}\tilde a_k Q_{\vec m_k}(x,y)
            $$
            where $\tilde a_{k+1}=\prodesc{xQ_{\vec n}(x,y)}{P_{\vec m_k}(x,y)}, \ \ k=|\vec n|-d_{\vec n},\dots,|\vec n| + (d_{\vec n}+1)r-1$.
            
            \item Consider $\{\vec m_k: k=|\vec n|-(d_{\vec n}+1),\dots,|\vec n| + (d_{\vec n}+2)r\}$ a path of neighbour multi-indices where $|\vec m_k| = k$, $\vec m_{|\vec n|}=\vec n$ and $\vec m_{|\vec n| + (d_{\vec n}+2)r}= (n_1+d_{\vec n}+2,\dots,n_r+d_{\vec n}+2)$. For every k, let $Q_{\vec m_k}$ be the bivariate Type I Function and $P_{\vec m_k}$ the bivariate Type II MOP.  Then:
            $$
                y Q_{\vec n}(x,y) = Q_{\vec m_{|\vec n|-(d_{\vec n}+1)}}(x,y) + \sum_{k=|\vec n|-d_{\vec n}}^{|\vec n|+(d_{\vec n}+2)r}\tilde b_k Q_{\vec m_k}(x,y)
            $$
            where $\tilde b_{k+1}=\prodesc{yQ_{\vec n}(x,y)}{P_{\vec m_k}(x,y)}, \ \ k=|\vec n|-d_{\vec n}-1,\dots,|\vec n| + (d_{\vec n}+2)r-1$.
        \end{enumerate}
    \end{theorem}
    \begin{proof}
        Again, only the first item will be proven, since the second one is analogous but having in mind the differences between the multiplication by `$x$' and `$y$' and employing the second equation of Lemma \ref{lemma:pi-xyP}.

        Let $\{\vec m_k: k=0,\dots,|\vec n| +  (d_{\vec n}+1)r\}$ be a path of neighbour multi-indices  where $|\vec m_k| = k$, $\vec m_0 =0$, $\vec m_{|\vec n|}=\vec n$ and $\vec m_{|\vec n|+(d_{\vec n}+1)r}=(n_1+d_{\vec n}+1,\dots,n_r+d_{\vec n}+1)$.

        We have $x Q_{\vec n}(x,y) = \sum_{j=0}^r x A_{\vec n,j}(x,y) w_j(x,y)$. Fixing $j\in\{1,\dots,r\}$, 
        $$\pi(\mdeg(x A_{\vec n,j}))=\pi(\mdeg(A_{\vec n,j}))+\deg(A_{\vec n,j})+1\leq n_j - 1 + d_{\vec n} + 1 = n_j + d_{\vec n}.$$

        In order to express $x A_{\vec n,j}$ as a linear combination of bivariate Type I MOP associated to other multi-indices, we need at least a multi-index whose $j$-th component is  $n_j + d_{\vec n} +1$, so that the $j$-th bivariate Type I MOP has an appropriate leading term. Applying this for every $j=1,\dots,r$, consider the multi-index $(n_1+d_{\vec n}+1,\dots,n_r+d_{\vec n}+1)$, whose norm is $|\vec n|+(d_{\vec n}+1)r$.
        So, we have 
        $$x Q_{\vec n} = \sum_{k=0}^{|\vec n|+(d_{\vec n}+1)r}\tilde a_k Q_{\vec m_k},$$
        since $Q_{\vec 0}=0$,
        $$
        x Q_{\vec n} = \sum_{k=1}^{|\vec n|+(d_{\vec n}+1)r}\tilde a_k Q_{\vec m_k}.
        $$
        Multiplying both sides by $P_{\vec m_l}$, we get
        $$
        \prodesc{x Q_{\vec n}}{P_{\vec m_l}} = \sum_{k=1}^{|\vec n|+(d_{\vec n}+1)r}\tilde a_k \prodesc{Q_{\vec m_k}}{P_{\vec m_l}},
        $$
        and applying the biorthogonality relation \eqref{eq:biorthogonality-2-app},
        $$
        \prodesc{x Q_{\vec n}}{P_{\vec m_l}} = \tilde a_{l+1}, \ \ \ l=0,\dots,|\vec n|+(d_{\vec n}+1)r-1.
        $$
        Finally, we will prove that $\tilde a_1= \cdots = \tilde a_{|\vec n|-d_{\vec n}-1}=0$, that is, $\tilde a_{k+1}=\prodesc{x Q_{\vec n}}{P_{\vec m_k}}=0$ for \linebreak $k=0,\dots,|\vec n|-d_{\vec n}-2$.

        Due to Type I orthogonality, $\prodesc{x Q_{\vec n}}{P_{\vec m_k}}=\prodesc{Q_{\vec n}}{x P_{\vec m_k}}=0$ if $\pi(\mdeg(x P_{\vec m_k}))\leq |\vec n|-2$. If $k\leq|\vec n|-d_{\vec n}-2$, then $\pi(\mdeg(P_{\vec m_k}))\leq |\vec n|-d_{\vec n}-2$ and $\deg(P_{\vec m_k})\leq d_{\vec n} - 1$. Applying Lemma \ref{lemma:pi-xyP}:
        $$
        \pi(\mdeg(x P_{\vec m_k})) = \pi(\mdeg(P_{\vec m_k}))+\deg(P_{\vec m_k})+1\leq  |\vec n|-d_{\vec n}-2 + d_{\vec n} - 1 +1 = |\vec n|-2.
        $$
        And we have proven $\tilde a_{k+1}=0$ for $k=0,\dots,|\vec n|-d_{\vec n}-2$. Lastly, if $k=|\vec n|-d_{\vec n}-1$:
        $$
        \pi(\mdeg(x P_{\vec m_k})) = |\vec n|-d_{\vec n}-1 + d_{\vec n} - 1 +1 = |\vec n|-1,
        $$
        so that, due to \eqref{eq:typeI-conditions-function-2-app}, $\tilde a_{|\vec n|-d_{\vec n}}=1$ and we finally get the relation.

    \end{proof}

    In the univariate case, it is possible to express the function $x Q_{\vec n}$ as a linear combination of Type I functions associated to the multi-index $\vec n$ itself, one multi-index of lower norm, and $r$ multi-indices of higher norm. In this case, following the ideas of Theorem \ref{th:NNR-xyPn}, $xQ_{\vec n}(x,y)$ can be written as a linear combination of Type I functions for $\vec n$, $d_{\vec n}$ multi-indices of lower norm, and $(d_{\vec n}+1)r$ multi-indices of higher norm.

    In case of multiplying by `$y$', we have chosen $(d_{\vec n}+2)r$ multi-indices of higher norm, while in \eqref{eq:NNR-xPn} we employed only $(d_{\vec n}+1)r$. Also, multiplying by `$x$' only $d_{\vec n}$ multi-indices of lower norm are used, whereas in the second item we needed one more multi-index of lower norm.

\section{Relation with the univariate case}\label{SecRelation}

In this section we will show how the product of two univariate MOP for different systems of measures behaves as a bivariate MOP regarding our definitions. Let us consider $(\mu_1,\dots,\mu_{r_1})$ and $(\phi_1,\dots,\phi_{r_2})$, two perfect systems of $r_1$ and $r_2$ univariate measures and two multi-indices $\vec n=(n_1,\dots,n_{r_1})\in\NN^{r_1}$, $\vec m=(m_1,\dots,m_{r_2})\in\NN^{r_2}$. We will take the Type II MOP $P_{\vec n}(x)$ with respect to the system $\mu_1,\dots,\mu_{r_1}$, and the Type II MOP $P_{\vec m}(y)$ with respect to the system $\phi_1,\dots,\phi_{r_2}$. Then, we have 
\begin{align}
    P_{\vec n}(x) &= x^{|\vec n|} + \mathcal{O} (x^{|\vec n|-1}) & P_{\vec m}(y) &= y^{|\vec m|} + \mathcal{O} (y^{|\vec m|-1}),
\end{align}
so that
$$
R(x,y) = P_{\vec n}(x) P_{\vec m}(y) = x^{|\vec n|} y^{|\vec m|} + \mathcal{O} (x^{|\vec n|-1}y^{|\vec m|-1}).
$$

Thus, $\mdeg(R(x,y))=(|\vec n|,|\vec m|)$ and $\deg(R(x,y))=|\vec n|+|\vec m|$. If we denote as $\prodesc{\cdot}{\cdot}_{i,j}$ the inner product with respect to the bivariate measure $\mu_i(x)\times \phi_j(y)$, $i=1,\dots,r_1, \ j=1,\dots,r_2$, then the following equation holds:
\begin{equation}
    \label{eq:product-of-univariate-mop}
    \begin{split}
        \prodesc{P_{\vec n}(x)P_{\vec m}(y)}{x^t y^s}_{i,j} &= \int_{\RR}\int_{\RR} P_{\vec n}(x)P_{\vec m}(y)x^t y^s d\mu_i(x) d\phi_j(y) \\
        &= \left(\int_{\RR}P_{\vec n}(x)x^t d\mu_i(x) \right)\left(\int_{\RR} P_{\vec m}(y) y^s  d\phi_j(y)\right) \\
        &= 0 \ \ \ 
            t = 0,\dots,n_i-1, \quad\text{ or } s = 0,\dots,m_j-1.
    \end{split}
\end{equation}

As $\mdeg(R(x,y))=(|\vec n|,|\vec m|)$, then $R(x,y)$ has $\pi(|\vec n|,|\vec m|)$ coefficients $c_0,...,c_{\pi(|\vec n|,|\vec m|)-1}$:
$$
R(x,y) = x^{|\vec n|} y^{|\vec m|} + c_{\pi(|\vec n|,|\vec m|)-1} x^{|\vec n|+1} y^{|\vec m|-1}  + \cdots + c_0.
$$
As a result, in order to get a bivariate multiple orthogonality relation, we need to find \linebreak $\vec v\in \NN^{r_1\times r_2}$ s.t. $|\vec v|=\pi(|\vec n|,|\vec m|)$ and satisfying $\prodesc{R(x,y)}{x^t y^s}_{i,j} = 0$ if $\pi(t,s)\leq v_{ij}-1,$ $i=1,\dots,r_1, \ j=1,\dots,r_2$.

On the other hand, due to \eqref{eq:product-of-univariate-mop}, observe that
$$
\prodesc{R(x,y)}{x^t y^s}_{i,j} = 0  \text{ if } \pi(t,s) \leq \pi(n_i,m_j)-1.
$$
Then, denoting $\tilde v = (\tilde v_{ij}=\pi(n_i,m_j), i=1,\dots,r_1, \ j=1,\dots,r_2)$, we know
$$
\prodesc{R(x,y)}{x^t y^s}_{i,j} = 0, \ \ \pi(t,s)\leq \tilde v_{ij}-1,\ \  i=1,\dots,r_1, \ j=1,\dots,r_2.
$$
As a consequence, it is essential for the orthogonality relation that $v_{ij}\leq \tilde v_{ij}$, otherwise the orthogonality relation would not hold. In this moment, the question is, is there any multi-index $\vec v=(v_{ij})\in \NN^{r_1\times r_2}$ so that $v_{ij}\leq \tilde v_{ij}$ and $|\vec v|=\pi(|\vec n|,|\vec m|)$?

We have proven an affirmative answer for this question for the particular case $r_1=r_2=2$, so that we consider the systems of univariate measures $(\mu_1(x),\mu_2(x))$ and $(\phi_1(y),\phi_2(y))$, together with the bivariate system $(\mu_1(x)\phi_1(y), \mu_1(x)\phi_2(y), \mu_2(x)\phi_1(y), \mu_2(x)\phi_2(y))$. From this point, we will remove the dependences for `$x$' and `$y$', assuming $\mu_i$ is associated to `$x$' and $\phi_j$ to `$y$'.

\begin{propo}
    Consider two systems of univariate measures $(\mu_1,\mu_2)$ and $(\phi_1,\phi_2)$, and the bivariate system composed of the product of the univariate ones $(\mu_1\phi_1, \mu_1\phi_2, \mu_2\phi_1, \mu_2\phi_2)$. Let $\vec n=(n_1,n_2),\vec m=(m_1,m_2)\in\NN^2$ be two normal multi-indices for the systems $(\mu_1,\mu_2)$ and $(\phi_1,\phi_2)$, respectively. Define $\tilde v = (\tilde v_{ij}) =(\pi(n_1,m_1),\pi(n_1,m_2),\pi(n_2,m_1),\pi(n_2,m_2))$. Then, it is always possible to find a multi-index $\vec v\in\NN^4$ such that $|\vec v|=\pi(|\vec n|,|\vec m|)$ and $\vec v \leq \tilde v$ componentwise.
\end{propo}
\begin{proof}
    We have to prove that
    \begin{equation}
        \label{eq:product-0}
        |\tilde v|=\pi(n_1,m_1)+\pi(n_1,m_2)+\pi(n_2,m_1)+\pi(n_2,m_2) \geq \pi(n_1+n_2,m_1+m_2)= \pi(|\vec n|,|\vec m|).
    \end{equation}
    If we prove this inequality, it is easy to find a multi-index $\vec v\in\NN^4$ so that $|\vec v|=\pi(|\vec n|,|\vec m|)$ by substracting the exceeding units from the components of $\tilde v$. Then, for all $n_1,n_2,m_1,m_2 \in \NN$:
    \begin{itemize}
        \item First step: Prove that 
        \begin{equation}
            \label{eq:product-1}
            \pi(n_1+n_2,m_1+m_2)=\pi(n_1,m_1)+\pi(n_2,m_2)+(n_1+m_1)(n_2+m_2).
        \end{equation}
        Applying the definition of the Cantor Pairing Function  \eqref{eq:CPF} this equality is easy to prove:
        \begin{equation*}
            \begin{array}{c}\vspace{5pt}
                \pi(n_1+n_2,m_1+m_2) = \frac 1 2 (n_1+n_2+m_1+m_2)(n_1+n_2+m_1+m_2+1)+m_1+m_2 \\ \vspace{5pt}
                = \frac{1}{2}((n_1+m_1)+(n_2+m_2))((n_1+m_1)+(n_2+m_2)+1)+m_1+m_2 \\ \vspace{5pt}
                = \frac{1}{2}((n_1+m_1)^2+(n_1+m_1)+ 2(n_1+m_1)(n_2+m_2)+(n_2+m_2)^2+(n_2+m_2)) + m_1+m_2 \\ \vspace{5pt}
                = \frac{1}{2}((n_1+m_1)(n_1+m_1+1)+ 2(n_1+m_1)(n_2+m_2)+(n_2+m_2)(n_2+m_2+1)) + m_1+m_2 \\ \vspace{5pt}
                = \frac{1}{2}(n_1+m_1)(n_1+m_1+1)+m_1 + \frac{1}{2}(n_2+m_2)(n_2+m_2+1)+m_2 + (n_1+m_1)(n_2+m_2) \\
                = \pi(n_1,m_1)+\pi(n_2,m_2)+(n_1+m_1)(n_2+m_2).
            \end{array}
        \end{equation*}
        \item Second step: Prove that  
        \begin{equation}
            \label{eq:product-2}
            (n_1+m_1)(n_2+m_2)\leq \pi(n_1,m_2)+\pi(n_2,m_1).
        \end{equation}
        \eqref{eq:product-2} is equivalent to 
        \begin{equation*}
            \begin{aligned}
                n_1 n_2 + n_1 m_2 +  n_2 m_1 + m_1 m_2 & \leq \\ \leq\frac 1 2(n_1+m_2)(n_1+m_2+1)&+\frac 1 2(n_2+m_1)(n_2+m_1+1)+m_1+m_2,
            \end{aligned}
        \end{equation*}
        which is in turn equivalent to
        \begin{equation*}
            \begin{aligned}
                2(n_1 n_2 + n_1 m_2 + n_2 m_1 + m_1 m_2) & \leq \\ \leq(n_1+m_2)(n_1+m_2+1)&+(n_2+m_1)(n_2+m_1+1)+2 m_1+2 m_2
            \end{aligned}
        \end{equation*}

        Expanding both members of the inequality, we get
        \begin{equation*}
            \begin{aligned}
        n_1^2+n_1+m_2^2+3m_2+n_2^2+n_2  +m_1^2+3 m_1 - 2 n_1 n_2& - 2 m_1 m_2 = \\
        (n_1-n_2)^2+(m_1-m_2)^2 & + n_1+n_2 +  3(m_1+m_2)\geq 0.
            \end{aligned}
        \end{equation*}
        So, inequality \eqref{eq:product-2} holds.
    \end{itemize}
    With equations \eqref{eq:product-1} and \eqref{eq:product-2} it is possible to prove \eqref{eq:product-0} easily:
    \begin{equation*}
        \begin{split}
            \pi(n_1+n_2,m_1+m_2)= &\pi(n_1,m_1)+\pi(n_2,m_2)+(n_1+m_1)(n_2+m_2) \ \ \ \text{(due to \eqref{eq:product-1})} \\
            \leq & \pi(n_1,m_1)+\pi(n_2,m_2) + \pi(n_1,m_2)+\pi(n_2,m_1)\ \ \ \text{(due to \eqref{eq:product-2}).} 
        \end{split}
    \end{equation*}

\end{proof}

 \begin{remark}
     The election of $\vec v$ is not unique, in fact, there exist many possibilities for the multi-index $\vec v$. However, the bivariate  Type II MOP $P_{\vec v}(x,y)$ is the same for every possible $\vec v$, and it is equal to $R(x,y)=P_{\vec n}(x)Q_{\vec m}(y)$. It also happens in the univariate case. For a fixed system of measures, a polynomial can be a Type II MOP for different multi-indices.
\end{remark}

In order to clarify this fact, we present an example.

\begin{ex}
    We choose Laguerre measures in order to create, firstly, univariate multiple Laguerre Polynomials of the First Kind \citep[Section 3.6.1]{VA20}. Then, we present the following Laguere measures
     \begin{align*}
         d\mu_1(x) &= x e^{-x} &  d\mu_2(x) &= x^{2.2} e^{-x}, \\
         d\phi_1(y) &= y^{2.3} e^{-y} &  d\phi_2(y) &= y^{3.4} e^{-y} .
     \end{align*}
    Let $\vec n = (0,1)$, then $P_{(0,1)}(x)=x-4.4$ is the Type II MOP with respect to the system $(\mu_1,\mu_2)$. Now, let us take $\vec m = (1,0)$, so that $P_{(1,0)}(y)=y-3.3$ is the Type II MOP with respect to the system $(\phi_1,\phi_2)$. The product polynomial is 
    $$R(x,y)=P_{(0,1)}(x)P_{(1,0)}(y)= xy -4.4y - 3.3 +14.52.$$
    Now, we consider the bivariate measures $(\mu_1\phi_1,\mu_1\phi_2,\mu_2\phi_1,\mu_2\phi_2)$.
    As $|\vec n|=|\vec m|=1$, $\pi(|\vec n|,|\vec m|)=\pi(1,1)=4$ and
    $$\tilde v = (\pi(0,1),\pi(0,0),\pi(1,1),\pi(1,0)) = (2,0,4,1).$$
    Then, we need $\vec v\in\NN^4$ s.t. $|\vec v|=4$ and $\vec v\leq \tilde v$. For example, $\vec v_1 = (1,0,2,1)$; and of course,
     $$
     P_{(1,0,2,1)}(x,y) = xy -4.4y - 3.3 +14.52.
     $$
     But, in fact, we can choose $\vec v_2 = (2,0,1,1)\neq \vec v_1$, however, the result is the same:
     $$
     P_{(2,0,1,1)}(x,y) = xy -4.4y - 3.3 +14.52.
     $$
 \end{ex}
 
Nevertheless, there is a problem with the valid multi-indices, because Proposition \ref{prop:equivalence-2-app} is not satisfied for every possible multi-index $\vec v\in\NN^{4}$. As a consequence, the system of measures $(\mu_1 \phi_1, \mu_1 \phi_2, \mu_2 \phi_1, \mu_2 \phi_2)$ and, more generally, the bivariate systems composed of products of two perfect systems of univariate measures are not perfect. See the example below.
 
\begin{ex}
If we take $\vec n = (1,1)$ and $\vec m = (1,1)$, we have $|\vec n|= |\vec m| = 2$ and $\pi(|\vec n|,|\vec m|)=\pi(2,2)=12$. In addition,
    $$
    \tilde v = (\pi(1,1),\pi(1,1),\pi(1,1),\pi(1,1)) = (4,4,4,4)
    $$
So, we can take, for example $\vec v = (3,3,3,3)$, but $\det M_{(3,3,3,3)}=0$, so $(3,3,3,3)$ is not a normal index in the product system of measures. Despite this, if we choose $\vec v = (4,3,3,2)$, then  $\det M_{(4,3,3,2)}\neq 0$.
\end{ex}

In fact, for many multi-indices $\vec n\in\NN^4$, the determinant of the matrix $M_{\vec n}$ can be expressed as a product of determinants of matrices associated to multi-indices in the univariate case. See the following example.

\begin{ex}
    We will denote the moments of the measures as
    \begin{align*}
        a_k &= \int_\RR x^k d\mu_1(x) & b_k &=\int_\RR y^k d\phi_1(y)  \\
        c_k &= \int_\RR x^k d\mu_2(x) & d_k &=\int_\RR y^k d\phi_2(y)  \\
    \end{align*}
    Then,
    $$
    \det(M_{(0,3,1,1)}) \propto  d_0 \begin{vmatrix}
        a_0 & a_1 & c_0 \\
        a_1 & a_2 & c_1 \\
        a_2 & a_3 & c_2 
        \end{vmatrix} \begin{vmatrix}
            a_0 & c_0 \\
            a_1 & c_1  
        \end{vmatrix} \begin{vmatrix}
            b_0 & d_0 \\
            b_1 & d_1  
        \end{vmatrix} \begin{vmatrix}
            d_0 & d_1 \\
            d_1 & d_0  
        \end{vmatrix},
    $$
    which is the product of the determinants of the matrices associated to the multi-indices $(2,1)$ and $(1,1)$ for the system $(\mu_1,\mu_2)$ and $(0,1)$, $(1,1)$ and $(0,2)$ for the system $(\phi_1,\phi_2)$.
    $$
    \det(M_{(2,0,1,1)}) \propto  c_0 b_0^2 \begin{vmatrix}
        a_0 & a_1 & c_0 \\
        a_1 & a_2 & c_1 \\
        a_2 & a_3 & c_2 
        \end{vmatrix} \begin{vmatrix}
            b_0 & d_0 \\
            b_1 & d_1  
        \end{vmatrix},
    $$
    which is the product of the determinants of the matrices associated to the multi-indices $(0,1)$ and $(2,1)$ for the system $(\mu_1,\mu_2)$ and $(1,0)$ (twice) and $(1,1)$ for the system $(\phi_1,\phi_2)$.

    This pattern continues for many other multi-indices.
\end{ex}
 
\section{Conclusions and state of art}\label{SecConc}

We have presented two definitions for Type I and Type II MOPs using bivariate polynomials and based in T. Koornwinder's perspective given in \citep{koornwinder}. Also, we have proven the usability of this definitions by extending some results like the existence and equivalence of MOP employing moment matrices, the biorthogonality relation and the Nearest Neighbour Recurrence Relations for bivariate Type I and Type II MOP, comparing every single result with its univariate analogous.

For the future, we would like to employ these definitions and their properties to extend some Multiple Orthogonal Polynomials applications to the bivariate case. For example, employ the extensions of Padé approximants to the multivariate case by A. A. M. Cuyt in \citep{Cuy83, Cuy86, Cuy99} and apply these tools to find simultaneous rational approximations for several bivariate functions. Additionally, we want to find some examples of perfect systems of measures (\citep{kounchev2011multidimensionalchebyshevsystems}), together with extensions of some known families of MOP like Jacobi-Piñeiro, Jacobi-Angelesco or Multiple Laguerre.

\section*{Acknowledgements}

This work has been partially supported by the grants ``PID2023-149117NB-I00'' and ``CEX 2020-001105-M'' funded by ``MCIN/AEI/10.13039/501100011033'', and Research Group Goya-384, Spain.

\bibliographystyle{plain}
\bibliography{references}{}

\begin{thebibliography}{10}

\bibitem{Aptekarev}
A.I. Aptekarev.
\newblock {Multiple orthogonal polynomials}.
\newblock {\em Journal of Computational and Applied Mathematics}, 99(1):423--447, 1998.
\newblock Proceeding of the VIIIth Symposium on Orthogonal Polynomials and Their Application.

\bibitem{CB00}
B.~Benouahmane and A.~A.~M. Cuyt.
\newblock {Multivariate orthogonal polynomials, homogeneous Pad{\'e} approximants and Gaussian cubature}.
\newblock {\em Numerical Algorithms}, 24(1):1--15, 2000.

\bibitem{BP17}
C.~F. Bracciali and T.~E. P\'erez.
\newblock Bivariate orthogonal polynomials, 2{D} {T}oda lattices and {L}ax-type pairs.
\newblock {\em Appl. Math. Comput.}, 309:142--155, 2017.

\bibitem{BP23}
C.~F. Bracciali and M.~A. Pi{\~n}ar.
\newblock {On multivariate orthogonal polynomials and elementary symmetric functions}.
\newblock {\em Numerical Algorithms}, 92(1):183--206, 2023.

\bibitem{Foulquie23}
A.~Branquinho, A.~Foulqui\'e-Moreno, and M.~Ma\~nas.
\newblock {Oscillatory banded Hessenberg matrices, multiple orthogonal polynomials and Markov chains}.
\newblock {\em Physica Scripta}, 98(10):105223, October 2023.

\bibitem{branquinho2021}
A.~Branquinho, A.~Foulqui\'e-Moreno, M.~Ma\~nas, C.~\'alvarez Fern\'andez, and J.~E. Fern\'andez-D\'iaz.
\newblock {Multiple Orthogonal Polynomials and Random Walks}, 2021.

\bibitem{branquinho2021-2}
A.~Branquinho, A.~Foulqui\'e-Moreno, M.~Ma\~nas, and J.~E. Fern\'andez-D\'iaz.
\newblock {Hypergeometric Multiple Orthogonal Polynomials and Random Walks}, 2021.

\bibitem{coussement-2003}
E.~Coussement and W.~Van~Assche.
\newblock {Multiple orthogonal polynomials associated with the modified Bessel functions of the first kind}.
\newblock {\em Constructive Approximation}, 19(2):237--263, 3 2001.

\bibitem{CF24}
R.~Cruz-Barroso and L.~Fern\'andez.
\newblock {Orthogonal Laurent Polynomials of Two Real Variables}.
\newblock {\em Studies in Applied Mathematics}, 10 2024.

\bibitem{Cuy83}
A.~A.~M. Cuyt.
\newblock {Multivariate Pad\'e-approximants}.
\newblock {\em Journal of Mathematical Analysis and Applications}, 96(1):283--293, 1983.

\bibitem{Cuy86}
A.~A.~M. Cuyt.
\newblock {Multivariate Pad\'e approximants revisited}.
\newblock {\em BIT Numerical Mathematics}, 26:71--79, 1986.

\bibitem{Cuy99}
A.~A.~M. Cuyt.
\newblock {How well can the concept of Pad\'e approximant be generalized to the multivariate case?}
\newblock {\em Journal of Computational and Applied Mathematics}, 105(1):25--50, 1999.

\bibitem{CLY16}
A.~A.~M. Cuyt, W.~Lee, and X.~Yang.
\newblock {On tensor decomposition, sparse interpolation and Pad{\'e} approximation}.
\newblock {\em Ja{\'e}n Journal of Approximation}, 8(1):33--58, 2016.

\bibitem{dunkl_xu_2014}
C.~F. Dunkl and Y.~Xu.
\newblock {\em {Orthogonal Polynomials of Several Variables}}.
\newblock Encyclopedia of Mathematics and its Applications. Cambridge University Press, 2 edition, 2014.

\bibitem{FI23}
L.~Fern\'andez and M.~D. de~la Iglesia.
\newblock {QBD Processes Associated with Jacobi--Koornwinder Bivariate Polynomials and Urn Models}.
\newblock {\em Mediterranean Journal of Mathematics}, 20(6):290, Aug 2023.

\bibitem{FI21}
L.~Fern\'andez and M.D. de~la Iglesia.
\newblock {Quasi-birth-and-death processes and multivariate orthogonal polynomials}.
\newblock {\em Journal of Mathematical Analysis and Applications}, 499(1):125029, 2021.

\bibitem{Ismail}
M.~E.~H. Ismail.
\newblock {\em Classical and quantum orthogonal polynomials in one variable}, volume~98 of {\em Encyclopedia of mathematics and its applications}.
\newblock Cambridge University Press, Cambridge (England), 2005.

\bibitem{koornwinder}
T.~Koornwinder.
\newblock {Two-Variable Analogues of the Classical Orthogonal Polynomials}.
\newblock In Richard~A. Askey, editor, {\em Theory and Application of Special Functions}, pages 435--495. Academic Press, 1975.

\bibitem{kounchev2011multidimensionalchebyshevsystems}
O.~Kounchev.
\newblock {Multidimensional Chebyshev Systems - just a definition}, 2011.

\bibitem{Kuijlaars_2008}
A.~B.~J. Kuijlaars, A.~Mart\'inez-Finkelshtein, and F.~Wielonsky.
\newblock {Non-Intersecting Squared Bessel Paths and Multiple Orthogonal Polynomials for Modified Bessel Weights}.
\newblock {\em Communications in Mathematical Physics}, 286(1):217--275, October 2008.

\bibitem{Lisi}
M.~Lisi.
\newblock Some remarks on the {C}antor pairing function.
\newblock {\em Matematiche (Catania)}, 62(1):55--65, 2007.

\bibitem{Andrei2022}
A.~Mart\'inez-Finkelshtein, R.~Orive, and J.~S\'anchez-Lara.
\newblock {Electrostatic Partners and Zeros of Orthogonal and Multiple Orthogonal Polynomials}.
\newblock {\em Constructive Approximation}, 58(2):271--342, December 2022.

\bibitem{martinezfinkelshtein-2016}
A.~Mart\'inez-Finkelshtein and W.~Van~Assche.
\newblock {WHAT IS...A multiple orthogonal polynomial?}
\newblock {\em Notices of the American Mathematical Society}, 63(09):1029--1031, 10 2016.

\bibitem{nikishin1991rational}
E.~M. Nikishin and V.~N. Sorokin.
\newblock {\em Rational approximations and orthogonality}, volume~92.
\newblock American Mathematical Society Providence, RI, 1991.

\bibitem{vanassche2006}
W.~Van~Assche.
\newblock {Pad\'e and Hermite-Pad\'e approximation and orthogonality}, 2006.

\bibitem{vanassche2011}
W.~Van~Assche.
\newblock {Nearest neighbor recurrence relations for multiple orthogonal polynomials}.
\newblock {\em Journal of Approximation Theory}, 163(10):1427--1448, 2011.

\bibitem{VA20}
W.~Van~Assche.
\newblock {\em {Orthogonal and Multiple Orthogonal Polynomials, Random Matrices, and Painlev\'e Equations}}, chapter 13 (Part II), pages 629--683.
\newblock Springer International Publishing, 2020.

\end{thebibliography}

\end{document}